\documentclass[12pt,a4paper]{amsart}
\usepackage{amsfonts, amsmath, amssymb, amscd, amsthm, graphicx}
\usepackage{epsfig, color}
\usepackage{hyperref} 

\makeatletter
\def\blfootnote{\xdef\@thefnmark{}\@footnotetext}
\makeatother

\newtheorem{thm}{Theorem}[section]
\newtheorem{cor}[thm]{Corollary}
\newtheorem{lem}[thm]{Lemma}
\newtheorem{prop}[thm]{Proposition}

\theoremstyle{definition}
\newtheorem{defn}[thm]{Definition}
\theoremstyle{remark}
\newtheorem{rem}[thm]{Remark}

\newfont{\eufm}{eufm10}

\newcommand{\G}{\Gamma (G, X\cup \mathcal H)}
\newcommand{\Gs}{\Gamma (G, S)}

\newcommand{\Hl}{\{ H_\lambda \} _{\lambda \in \Lambda } }

\newcommand{\e}{\varepsilon }
\newcommand{\N}{\mbox{\eufm N}}

\renewcommand{\phi}{{\rm Lab }}
\renewcommand{\kappa}{\varkappa}

\renewcommand{\ll}{\left\langle\hspace{-.7mm}\left\langle }
\newcommand{\rr}{\right\rangle\hspace{-.7mm}\right\rangle }

\newcommand{\cH}{\mathcal{H}}
\newcommand{\cW}{\mathcal{W}}
\renewcommand{\L}{{\rm L}}
\newcommand{\NN}{{\mathbb N}}
\newcommand{\Z}{{\mathbb Z}}
\newcommand{\apg}{\stackrel{G}{\approx}}
\newcommand{\napg}{\stackrel{G}{\not\approx}}

\newcommand{\CH}{C_H(E_G(H))}

\begin{document}

\title[Normal automorphisms of relatively hyperbolic groups]{
Normal automorphisms of relatively hyperbolic groups}

\author{A. Minasyan}
\address[Ashot Minasyan]{School of Mathematics,
University of Southampton, Highfield, Southampton, SO17 1BJ, United
Kingdom.}  \email{aminasyan@gmail.com}

\author{D. Osin }
\address[Denis Osin]{Department of Mathematics, Vanderbilt University, Nashville TN 37240, USA.}
\email{denis.osin@gmail.com}
\dedicatory{Dedicated to Professor A.L. Shmelkin on the occasion of his 70th birthday.}
\thanks{The first author was supported by the Swiss National Science
Foundation grant PP002-116899. 
The second author was supported by the
NSF grant DMS-0605093 and by the RFBR grant 05-01-00895.}

\begin{abstract}
  An automorphism $\alpha$ of a group $G$ is normal if it fixes every normal
subgroup of $G$ setwise. We give an algebraic description of normal
automorphisms of relatively hyperbolic groups. In particular, we prove that for
any relatively hyperbolic group $G$, $Inn(G)$ has finite index in the subgroup $Aut_n(G)$ of normal automorphisms. If, in addition, $G$ is non-elementary and has no non-trivial finite normal subgroups, then $Aut_n(G)=Inn(G)$. As an application, we show that $Out(G)$ is residually finite for every finitely generated residually finite group $G$ with more than one end.
\end{abstract}


\maketitle
\tableofcontents

%


\section{Introduction}


Recall that an automorphism $\alpha \in Aut (G)$ of a group $G$ is
said to be {\it normal} if $\alpha(N)=N$ for every normal subgroup  $N$ of $G$.
The subset of normal automorphisms, denoted by $Aut_n(G)$,
is clearly a subgroup of $Aut (G)$. Obviously every inner
automorphism is normal. Throughout this paper we denote by
$Out_n(G)$ the quotient group $Aut_n(G)/Inn(G)$.

The study of normal automorphisms originates from the result of
Lubotzky stating that $Out_n(G)$ is trivial for any non-abelian free
group \cite{Lub}. Since then similar results have been proved for many other
classes of groups. For example, $Out_n(G)=\{ 1\} $ for non-trivial
free products \cite{Nesh}, fundamental groups of closed surfaces of
negative Euler characteristic \cite{BKZ}, non-abelian free Burnside
groups of large odd exponent \cite{Cher}, non-abelian free solvable
groups \cite{Rom}, and free nilpotent group of class $c=2$ (for
$c\ge 3$ this is not true) \cite{End}.  On the other hand, every
group can be realized as $Out (G)$ for a suitable simple group
$G$ \cite{DGG}. Since every automorphism of a simple group is normal,
every group appears as $Out_n(G)$ for some $G$. Furthermore, every countable group can be realized as
$Out_n(G)$ for some finitely generated group $G$ \cite{CC,Obr}.

The main goal of this paper is to study normal automorphisms of
relatively hyperbolic groups. The notion of a relatively
hyperbolic group was originally suggested by Gromov \cite{Gro} and
has been elaborated in many papers since then
\cite{Bow,DS,F,Hru,RHG,Yam}. The class of relatively hyperbolic groups
includes (ordinary) hyperbolic groups, fundamental groups of
finite-volume complete Riemannian manifolds of pinched negative
curvature \cite{Bow,F}, groups acting freely on $\mathbb R^n$-trees
\cite{Gui} (in particular, limit groups arising in the solutions of
the Tarski problem \cite{KM,Sel}), non-trivial free products and their small
cancellation quotients \cite{RHG}, groups acting geometrically on $CAT(0)$
spaces with isolated flats \cite{HK}, and many other examples.

In this paper we neither assume relatively hyperbolic groups to be finitely generated nor
the collection of peripheral subgroups to be finite. (The reader is referred to the next section for the precise definition.) However we do assume that all peripheral subgroups are proper to exclude the case of a group hyperbolic relative to itself. Further on, we will say that a group $G$ {\it non-elementary}, if it is not virtually cyclic.

In general $Out_n(G)$ is not necessarily trivial even
for ordinary hyperbolic groups.  Indeed, it is known (see \cite{Sah}) that certain finite groups
$L$ possess non-inner automorphisms which map every element to its conjugate. 
One can therefore
construct many hyperbolic groups $G$ with non-trivial $Out_n(G)$ by
taking any hyperbolic group $H$ and considering the direct product
$G=H\times L$. The first result of our paper shows that non-trivial finite
normal subgroups are essentially the only sources
of non-inner normal automorphisms.

More precisely, every relatively hyperbolic group $G$
contains a unique maximal finite normal subgroup
(see Corollary \ref{KG}). We denote it by $E(G)$. Further let $C(G)$ denote the centralizer of $E(G)$ in $G$.

\begin{thm}\label{main}
Suppose that $G$ is a non-elementary relatively hyperbolic group
and $\alpha \in Aut_n(G)$. Then there exist an element $w\in G$ and a set map
$\e\colon G\to E(G)$ such that $\e (C(G))=\{ 1\}$ and
$\alpha (g)=wg\e(g)w^{-1}$ for every $g\in G$.
\end{thm}

In fact, Theorem \ref{main} is a particular case of a more general
result about normal automorphisms of subgroups of relatively
hyperbolic groups (see Theorem \ref{NormAutSubgr}).
The corollary below follows easily from Theorem \ref{main} and the observation that
$C(G)$ has finite index in $G$ being the centralizer of a
finite normal subgroup.

\begin{cor}\label{cor1}
Suppose that $G$ is a  relatively hyperbolic group. Then the
following hold.
\begin{enumerate}
\item[(a)] $Out_n(G)$ is finite.

\item[(b)] If $G$ is non-cyclic and contains no non-trivial finite normal subgroups, then $Out_n(G)=\{ 1\} $.
\end{enumerate}
\end{cor}

This corollary generalizes the results about free groups \cite{Lub}, free products \cite{Nesh}, and
surface groups \cite{BKZ} cited above. It also implies the result of Metaftsis and Sykiotis \cite{MS} stating
that for every non-elementary finitely generated relatively hyperbolic group
$G$, $Inn(G)$ has finite index in the group $Aut_c(G)$ of pointwise inner automorphisms of $G$.  Recall that an automorphism of $G$ is {\it pointwise inner}, if it preserves conjugacy classes. Clearly $Aut_c(G)\le Aut_n(G)$. Thus finiteness of $Out_n(G)$
implies that of $Aut_c(G)/Inn(G)$. The converse is not true in general. For instance, if $G$ is free
nilpotent of class at least $3$, we have $Aut_c(G) = Inn(G) $ while
$\left| Out_n(G) \right| = \infty$ \cite{End}.

It is also worth noting that our methods are quite different from those of \cite{MS}. Indeed
we use the group-theoretic version of Dehn surgery introduced in
\cite{GM1,GM2,CEP} and `component analysis' developed in
\cite{RHG,CC}, while Metaftsis and Sykiotis employed the Bestvina-Paulin
approach \cite{Best,Pau} based on ultralimits and group actions on
$\mathbb R$-trees.

In order to prove Theorem \ref{main}, we introduce a new subclass of automorphisms of any given group, and investigate it in the case of relatively hyperbolic groups.
\begin{defn}
Let $G$ be a group. We say that an automorphism $\varphi \in Aut(G)$ is {\it commensurating} if for every $g\in G$ there exist
$h \in G$ and $m,n \in \Z\setminus \{0\}$ such that $(\varphi (g))^n = h g^m h^{-1}$. In other words, $\varphi$ is commensurating if
for each $g \in G$, $\varphi(g)$ is
{\it commensurable} to $g$ in $G$ (see Definition \ref{def:commensurability}).
\end{defn}

It is clear that the subset $Aut_{comm}(G)$ of commensurating automorphisms of $G$  forms a subgroup of $Aut(G)$ and
$Inn(G)\le Aut_c(G) \le Aut_{comm}(G)$.

In Section \ref{sec:comm-aut} we study commensurating automorphisms of relatively hyperbolic groups and obtain
a complete description of them:

\begin{cor} \label{cor:descr_comm_aut}
Let $G$ be a non-elementary relatively hyperbolic group and $\varphi \in Aut(G)$. The following conditions are equivalent:
\begin{enumerate}
\item[(i)] $\varphi $ is commensurating;

\item[(ii)] there is a set map $\e : G \to E(G)$, whose restriction to $C(G)$ is a homomorphism,
and an element $w \in G$ such that for every $g \in G$, $\varphi(g)=w \left( g \e (g) \right)w^{-1}$.
\end{enumerate}

In particular, if $E(G)=\{ 1\}$, then every commensurating automorphism of $G$ is inner.
\end{cor}

In Section \ref{sec:Dehn_surgery}, using the algebraic version of Dehn filling, we show that each normal
automorphism of a relatively hyperbolic group must be commensurating. After this, Theorem \ref{main} follows quite quickly from
the above description of commensurating automorphisms.

Our methods can also be used to prove residual finiteness of some outer automorphism groups.
A well-known theorem of Baumslag
states that the automorphism group of a
finitely generated residually finite group is residually finite \cite{Bau}. In
general, the analogous property does not hold for the group of outer
automorphisms. Indeed, Bumagina and Wise
showed that every finitely presented group is realized as $Out (G)$ for
a suitable finitely generated residually finite group $G$ \cite{BW}.
However we prove that Baumslag's theorem does
have an `outer' analogue for groups with more than one end.
We refer to \cite{Stall71} for the geometric definition of ends, and
recall that the number of ends of a finitely generated group 
can be either $0$, $1$, $2$ or infinity. 

\begin{thm}\label{infends}
Let $G$ be a finitely generated residually finite group with more than one end.
Then $Out (G)$ is residually finite.
\end{thm}

An infinite finitely generated group $G$ has two ends if and only if it is virtually cyclic; 
and $G$ has infinitely many ends if and only if it splits non-trivially as
an amalgamated free product $A\ast _S B$ or an $HNN$-extension $A\ast_S$
over a finite group $S$ \cite{Stall71,Stall68}. 


Note that the condition demanding residual finiteness of $G$ in Theorem \ref{infends} cannot be removed.
Indeed, if $H$ is any finitely
generated group that has trivial center and is not residually finite, then the group $G=H * \Z$ has infinitely many ends and
$H$ is embedded into $Out(G)$ ($H$ acts on itself by conjugation and trivially on $\Z$, which gives rise to an action of $H$
by outer automorphisms on the free product $H * \Z=G$). Thus $Out(G)$ is not residually finite.

The standard way of proving residual finiteness of $Out(G)$ is based on
the following result of Grossman \cite{Gros}: if  a group $G$ is finitely generated and conjugacy separable, then $Aut(G)/Aut_c(G)$ is residually finite. In particular,
$Out (G)$ is residually finite whenever $G$ is finitely generated,
conjugacy separable, and $Aut_c(G)=Inn(G)$. Recall that a group $G$ is said to be {\it conjugacy
separable} if for any two non-conjugate elements $g,h\in G$ there
exists a homomorphism $\varphi\colon G\to K$, where $K$ is finite, such
that $\phi (g)$ and $\varphi(h)$ are not conjugate in $K$.

This approach has been successfully used to prove residual finiteness
of $Out(G)$, where $G$ is a free group of finite rank \cite{Gros},
the fundamental group of a closed surface \cite{Gros,AKT01}, the
fundamental group of a Seifert manifold with non-trivial boundary \cite{AKT03}, etc.
If $G$ is a finitely generated conjugacy separable non-elementary relatively
hyperbolic group, the above mentioned result from \cite{MS} implies that every
virtually torsion-free subgroup of $Out(G)$ is residually finite \cite[Theorem 1.1]{MS}.

However there is no hope to use Grossman's idea to prove Theorem \ref{infends}
since we only assume the group $G$ to be residually finite, which is much
weaker than conjugacy separability. Indeed there are many examples of finitely
generated residually finite groups that are not conjugacy separable (e.g., the group of unimodular matrices $GL(n,\mathbb{Z})$
for $n\ge 3$, see \cite{Rem}). To construct such an example with infinitely many ends,
we can simply take $G=H\ast \mathbb Z$, where $H$ is finitely generated, residually finite,
but not conjugacy separable. It is easy to show that $G$ will also be finitely generated, residually finite,
but not conjugacy separable.

Our approach is different and is based on the  following observation.
Let $Aut_n^f(G)$ denote the group of automorphisms of $G$ that stabilize
every normal subgroup of finite index (setwise). Then $Aut(G)/Aut_n^f(G)$
is residually finite for every finitely
generated group $G$. The following result plays the crucial role in the
proof of Theorem \ref{infends}. It also seems to be of independent interest.
Its proof essentially uses the fact that free products are hyperbolic relative
to their free factors, which allows us to employ the techniques developed in the proof of Theorem \ref{main}.

\begin{thm}\label{fp}
Suppose that $G=A\ast B$, where $A,B$ are non-trivial residually finite groups. Then $Aut_n^f(G)=Inn(G)$.
\end{thm}

\noindent {\bf Acknowledgments.} We are grateful to A. Klyachko and V. Yedynak
for useful discussions, and to the anonymous referee for his comments.


\section{Preliminaries}\label{sec:prelim}


\noindent{\bf Notation.} Given a group $G$ generated by a subset $S\subseteq G$, we denote by $\Gs $ the
Cayley graph of $G$ with respect to $S$ and by $|g|_S$ the word length of an element $g\in G$. If $p$
is a (combinatorial) path in $\Gs$, $\phi (p)$ denotes its label, $\L(p)$ denotes its length, $p_-$ and $p_+$ denote
its starting and ending vertex. The notation $p^{-1}$ will be used for the path in $\Gs$ obtained by traversing $p$ backwards.
By saying that $o=p_1\dots p_k$ is a cycle in $\Gs$ we will mean that $o$
is obtained as a consecutive concatenation of paths $p_1,\dots p_k$
such that $(p_{i+1})_-=(p_i)_+$ for $i=1,\dots,k-1$ and $(p_k)_+=(p_1)_-$.

For a word $W$ written in the alphabet $S$, $\|W\|$ will denote its length. For two words $U$ and $V$ we shall write $U \equiv V$ to denote
the letter-by-letter equality between them.
The normal closure of a subset $K\subseteq G$ in a group $G$ (i.e., the minimal normal subgroup of $G$ containing $K$)
is denoted by $\ll K\rr^G$, or simply by $\ll K\rr$ if omitting $G$ does not lead to a confusion.  For any group elements $g$ and $t$, $g^t$ denotes $t^{-1}gt$. If $A\subseteq G$ then $A^t=\{a^t~|~a\in A\}$. For a subgroup $H\le G$, $N_G(H)$ denotes the normalizer of a $H$ in $G$. That is, $N_G(H)=\{g \in G~|~gHg^{-1}=H\}$. Similarly by $C_G(H)$ we denote the centralizer of $H$ in $G$,
that is, $$ C_G(H)=\{ g\in G~|~ gh=hg, \; \forall\, h\in H\} .$$ Finally for two subsets $A,B$ of $G$, their product is the subset
$AB =\{ab~|~a \in A,b \in B\}$.

\paragraph{\bf Relatively hyperbolic groups.} In this paper we use
the notion of relative hyperbolicity which is sometimes called
strong relative hyperbolicity and goes back to Gromov \cite{Gro}.
There are many equivalent definitions of (strongly) relatively
hyperbolic groups \cite{Bow,DS,F,RHG}. We recall the isoperimetric
characterization suggested in \cite{RHG}, which is most suitable
for our purposes. That relative hyperbolicity in the sense of
\cite{Bow,F,Gro} implies relative hyperbolicity in the sense of
Definition \ref{def:rel_hyp_gp} stated below is essentially due to
Rebbechi \cite{Reb}. (Indeed it was proved in \cite{Reb} under the
additional technical condition that the groups under consideration
are finitely presented.) In the full generality this implication
and the converse one were proved in \cite{RHG}.

Let $G$ be a group, $\Hl $ -- a collection of {\it proper} subgroups of $G$, $X$ -- a
subset of $G$. We say that $X$ is a {\it relative generating set of
$G$ with respect to $\Hl $} if $G$ is generated by $X$ together with
the union of all $H_\lambda $. (In what follows we always assume $X$
to be symmetric.) In this situation the group $G$ can be regarded as
a quotient group of the free product
\begin{equation}
F=\left( \ast _{\lambda\in \Lambda } H_\lambda  \right) \ast F(X),
\label{F}
\end{equation}
where $F(X)$ is the free group with the basis $X$. If the kernel of
the natural homomorphism $F\to G$ is the normal closure of a subset
$\mathcal R$ in the group $F$, we say that $G$ has {\it relative
presentation}
\begin{equation}\label{G}
\langle X,\; H_\lambda, \lambda\in \Lambda \; \mid \; \mathcal R
\rangle .
\end{equation}
If $|X|<\infty $ and $|\mathcal R|<\infty $, the
relative presentation (\ref{G}) is said to be {\it finite} and the
group $G$ is said to be {\it finitely presented relative to the
collection of subgroups $\Hl $.}

Set
\begin{equation}\label{H}
\mathcal H=\bigsqcup\limits_{\lambda\in \Lambda} (H_\lambda
\setminus \{ 1\} ) .
\end{equation}
Given a word $W$ in the alphabet $X\cup \mathcal H$ such that $W$
represents $1$ in $G$, there exists an expression
\begin{equation}
W\stackrel{F}{=} \prod\limits_{i=1}^k f_i^{-1}R_i^{\pm 1}f_i \label{prod}
\end{equation}
with the equality in the group $F$, where $R_i\in \mathcal R$ and
$f_i\in F $ for $i=1, \ldots , k$. The smallest possible number
$k$ in a representation of the form (\ref{prod}) is called the
{\it relative area} of $W$ and is denoted by $Area^{rel}(W)$.

\begin{defn}[\cite{RHG}] \label{def:rel_hyp_gp}
A group $G$ is {\it hyperbolic relative to a collection of proper
subgroups} $\Hl $ if $G$ is finitely presented relative to $\Hl $
and there is a constant $C>0$ such that for any word $W$ in $X\cup
\mathcal H$ representing the identity in $G$, we have
\begin{equation}\label{isop}
Area^{rel}
(W)\le C\| W\| .
\end{equation}
The constant $C$ in (\ref{isop}) is called an {\it isoperimetric
constant} of the relative presentation (\ref{G}) and $\Hl $ is called the collection of {\it peripheral (or parabolic) subgroups} of $G$.
In particular, $G$ is an ordinary (Gromov) {\it hyperbolic group} if $G$ is hyperbolic relative to the trivial subgroup.
Later on by saying that a group $G$ is {\it relatively hyperbolic}, we will mean that there exists a collection of proper subgroups
$\{H_\lambda \le G~|~\lambda \in \Lambda\}$ such that $G$ is hyperbolic relative to $\Hl$.
\end{defn}

This definition is independent of the choice of the finite
generating set $X$ and the finite set $\mathcal R$ in (\ref{G})
(see \cite{RHG}).

\begin{lem}[\cite{RHG}, Thm. 1.4]\label{maln}
Let $G$ be a group hyperbolic relative to a collection of subgroups $\Hl $. Then the following conditions hold.
\begin{enumerate}
\item For every $\lambda, \mu \in \Lambda $, $\lambda \ne \mu $, and every $g\in G$, we have $|H_\lambda \cap H_\mu^g |<\infty $.

\item For every $\lambda \in \Lambda $ and $g\in G\setminus H_\lambda $, we have $|H_\lambda \cap H_\lambda ^g|<\infty $.
\end{enumerate}
\end{lem}

\paragraph{\bf Components.} Let $G$ be a group hyperbolic relative to a family of proper subgroups $\Hl$.
We recall some auxiliary terminology introduced in \cite{RHG}, which
plays an important role in our paper.

\begin{defn}
Let $q$ be a path in the Cayley graph $\G $. A (non-trivial)
subpath $p$ of $q$ is called an {\it $H_\lambda $-component} (or simply a {\it component}), if the label of $p$ is a word in the
alphabet $H_\lambda\setminus \{ 1\} $ and $p$ is not contained in a longer
subpath of $q$ with this property.  Two $H_\lambda $-components $p_1, p_2$
of a path $q$ in $\G $ are called {\it connected} if there exists a
path $c$ in $\G $ that connects some vertex of $p_1$ to some vertex
of $p_2$, and ${\phi (c)}$ is a word consisting of letters from
$H_\lambda\setminus\{ 1\} $. In algebraic terms this means that all
vertices of $p_1$ and $p_2$ belong to the same coset $gH_\lambda $
for a certain $g\in G$. Note that we can always assume that $c$ has
length at most $1$, as every non-trivial element of $H_\lambda
\setminus\{ 1\} $ is included in the set of generators.
\end{defn}

\paragraph{\bf Loxodromic elements and elementary subgroups.}
Recall that an element $g\in G$ is called {\it parabolic} if it is conjugate to an element of one of the
subgroups $H_\lambda $, $\lambda \in \Lambda $. An element is said to be {\it loxodromic} if it is not
parabolic and has infinite order. If $H$ is a subgroup of $G$, by $H^0 \subset H$ we will denote the set of all elements of
$H$ that are loxodromic in $G$.

Recall also that a group is {\it elementary} if it contains a cyclic
subgroup of finite index. The next result was obtained in \cite{ESBG}. The first part of the lemma is
well known in the context of convergence groups \cite{Tuk}. In particular, it follows from \cite{Tuk} 
and \cite{Yam} in case $G$ is finitely generated. (The latter assumption is only essential 
for \cite{Yam}.) 

\begin{lem}\label{Eg}
Suppose a group $G$ is hyperbolic relative to a collection of subgroups $\Hl$. Let $g$ be a
loxodromic element of $G$. Then the following conditions hold.
\begin{enumerate}
\item[(a)] There is a unique maximal elementary
subgroup $E_G(g)\le G$ containing $g$.

\item[(b)] $E_G(g)=\{ h\in G\mid \exists\, m\in \mathbb{N}~ \mbox{such that}~ h^{-1}g^mh=g^{\pm m}\} $.

\item[(c)] The group $G$ is hyperbolic relative to the collection
$\Hl\cup \{ E_G(g)\} $.
\end{enumerate}
\end{lem}

For finitely generated relatively hyperbolic groups, a lemma similar to Lemma \ref{Eg} (c) was also stated in \cite{Dah}. Namely it was claimed that if $G$ is a (finitely generated) relatively hyperbolic group and $Z$ is an infinite cyclic subgroup of $G$ such that $Z$ coincides with its normalizer, then $Z$ can be joined to the collection of peripheral subgroups of $G$ \cite[Lemma 4.4]{Dah}. We note that this is wrong even in case $G$ is an ordinary hyperbolic group.

The simplest counterexample is given by the group $$G=\langle z,c\mid c^3=1,\; zcz^{-1}=c^2\rangle $$ and the subgroup $Z=\langle z\rangle $. Obviously $G$ splits as $1\to C\to G\to \mathbb Z\to 1$, where $C=\langle c\rangle \cong \mathbb Z/3\mathbb Z$. In particular $G$ is hyperbolic (or, equivalently, hyperbolic relative to the trivial subgroup). It is straightforward to check that $Z$ coincides with its own normalizer in $G$. Indeed every element $g\in G$ has the form $z^kc^m$, where $k\in \mathbb Z$ and $m\in \{ 0,1,2\}$. If $m=1$, we have $$g^{-1}zg =(c^{-1}z^{-k})z(z^{k}c)=c^{-1}zc=c^{-1}(zcz^{-1})z=c^{-1}c^2z=cz\notin Z.$$ Similarly $g^{-1}zg\notin Z$ if $m=2$. On the other hand, $G$ is not hyperbolic relative to $Z$. Indeed $c^{-1}z^2c=z^2$ and hence $Z\cap c^{-1}Zc$ is infinite. This contradicts part (b) of Lemma \ref{maln}. Similarly for every (finitely generated) group $H$, the free product $G\ast H$ is hyperbolic relative to $H$, and the subgroup $Z$ provides a counterexample. Note that $E_G(z)=E_{G\ast H}(z)=G$, so applying Lemma \ref{Eg} (c) yields the correct result.

\paragraph{\bf Finite normal subgroups}
The following result was proved in \cite[Lemma 3.3]{AMO}.

\begin{lem}\label{EH} Let $H$ be a non-elementary subgroup of a relatively hyperbolic
group $G$. Suppose that $H^0\ne \emptyset $.  Then the subgroup $\displaystyle E_G(H)=\bigcap_{h \in
H^0} E_G(h)$ is the (unique) maximal finite subgroup of $G$
normalized by $H$.
\end{lem}

\begin{cor}\label{KG}
Let $G$ be a relatively hyperbolic group. Then $G$ possesses a unique maximal finite normal subgroup $E(G)$.
\end{cor}

\begin{proof}
If $G$ is finite then the statement is trivial. If $G$ contains an infinite normal cyclic subgroup $C$ of finite index, then
denote by $K$ the union of all finite normal subgroups of $G$. It is easy to see that $K$ is a torsion normal subgroup of $G$
(because a product of two finite normal subgroups is itself a finite normal subgroup). Since $K \cap C=\{1\}$, $K$
injects into the finite quotient $G/C$, hence $K$ is finite.

Finally, if $G$ is non-elementary, then
$G^0\ne \emptyset$ by \cite[Cor. 4.5]{ESBG}
(if $G$ is finitely generated, this also follows from \cite{Tuk} and \cite{Yam}).
It remains to apply Lemma \ref{EH} to the case $G=H$.
\end{proof}


\section{Special elements in relatively hyperbolic groups}


Let $G$ be a relatively hyperbolic group and let $H$ be a non-elementary subgroup of $G$ containing at least one loxodromic element.

\begin{defn} We say that an element $h\in H$ is $H$-{\it special} if $h$ is loxodromic in $G$
and $E_G(h)=\langle h\rangle\times E_G(H)$. The set of all $H$-special elements will be denoted by $S_G(H)$.
\end{defn}

Note that, by definition, any $g \in S_G(H)$ belongs to the centralizer $\CH$.
The result below was obtained in \cite[Lemma 3.8]{AMO}.

\begin{lem}\label{SG} If $G$ is a relatively hyperbolic group and $H \le G$ is a non-elementary subgroup such that $H^0\ne \emptyset$,
then $S_G(H)$ is non-empty.
\end{lem}

Special elements play a significant role in our approach to study automorphisms of relatively hyperbolic groups.
They represent a useful tool that helps to deal with the technical problems which may arise when the group under consideration
contains torsion.
The main goal of this section is to prove the following important statement:

\begin{prop}\label{fi}
Suppose that $G$ is a relatively hyperbolic group and $H \le G$ is a
non-elementary subgroup with  $H^0\ne \emptyset$. Then $C_H(E_G(H))$ is generated by the set
$S_G(H)$. In particular, $\langle S_G(H)\rangle $ has finite index
in $H$.
\end{prop}

Observe that the statement after `in particular' follows from the fact that the centralizer of a
finite subgroup of $G$, normalized by $H$, necessarily has finite index in $H$.

We begin with some auxiliary results. Let $G$ be a group hyperbolic relative to a family of proper
subgroups $\Hl$. If $G$ is infinite, it always contains a loxodromic element
\cite[Corollary 4.5]{ESBG}. The next lemma  provides us with a tool
for constructing such elements. It was proved in \cite[Lemma
4.4]{ESBG}.

\begin{lem}\label{ah}
Let $G$ be a group hyperbolic relative to a collection of subgroups $\Hl $. For any $\lambda \in \Lambda $
and $a\in G\setminus H_\lambda $, there exists a finite subset $\mathcal F\subseteq H_\lambda $ such that
if $h\in H_\lambda \setminus \mathcal F$, then $ah$ is loxodromic.
\end{lem}

Suppose that  $\Xi$ is a finite subset of $G$. Define $\cW (\Xi)$ to be the set of all
words $W$ over the alphabet $X \cup \cH$ that have the following form:
$$W \equiv x_0h_0x_1h_1 \dots x_l h_l x_{l+1},$$ where $ l \in \Z$, $l \ge -2$
(if $l=-2$ then $W$ is the empty word; if $l=-1$ then $W \equiv x_0$),
$h_i$ and $x_i$ are considered as single letters and
\begin{itemize}
\item[1)] $x_i \in X \cup \{1\}$, $i=0,\dots,l+1$,
and for each $i=0,\dots,l$, there exists $\lambda(i) \in \Lambda$ such that $h_i \in H_{\lambda(i)}$;
\item[2)] if $\lambda(i)=\lambda(i+1)$ then $x_{i+1} \notin H_{\lambda(i)}$ for each $i=0,\dots,l-1$;
\item[3)] $h_i \notin \Xi$, $i=0,\dots,l$.
\end{itemize}

The statement below was proved in \cite[Lemmas 6.3, 6.5]{CC}.
\begin{lem}\label{lem:conseq-conn} There is a finite subset $\Xi$ of $G$ such that the following holds.
Suppose that $o=rqr'q'$ is a cycle in $\Gamma(G,X \cup\cH)$ with $\phi(q),\phi(q') \in \cW(\Xi)$.
Set $C=\max\{ \L(r),\L(r')\}$.
\begin{itemize}
\item[(a)] If $l$ is the number of components of $q$, then at least $(l-6C)$ of
components of $q$ are connected to components of $q'$; and two distinct components of $q$ cannot be connected to the same
component of $q'$. Similarly for $q'$.

\item[(b)] For any $d\in \NN$ there exists
a constant $L=L(C,d) \in \NN$ such that if $\L(q)\ge L$ then there are $d$ consecutive components $p_s,\dots,p_{s+d-1}$ of $q$ and
$p'_{s'},\dots,p'_{s'+d-1}$ of $q'^{-1}$, so that $p_{s+i}$ is connected to $p'_{s'+i}$ for each $i=0,\dots,d-1$.
\end{itemize}
\end{lem}

Proposition \ref{fi} is an easy consequence of Lemma \ref{lem:spec-mod} below.
In the case when $G$ is an ordinary word hyperbolic group it was proved in \cite[Lemma 4.3]{paper3}.

\begin{lem}\label{lem:spec-mod} Suppose that $g \in S_G(H)$ and $x \in C_H(E_G(H)) \setminus E_G(g)$. Then there exists
$N_1 \in \NN$ such that $g^n x \in S_G(H)$ for any $n \in \Z$ with $|n| \ge N_1$.
\end{lem}

\begin{proof} By part (3) of Lemma \ref{Eg}, $G$ is hyperbolic relative to the collection of subgroups
$\Hl \cup \{E_G(g)\}$. Denote $\mathcal{H}'=\left(\cup_{\lambda\in \Lambda} H_\lambda \cup E_G(g)\right)\setminus \{1\} \subset G$.
After adding $x$ and $x^{-1}$ to the finite relative generating set of $G$, if necessary,
we can assume that $x^{\pm 1} \in X$.
Let $\mathcal F$ and $\Xi$ be the finite subsets of $G$ given by Lemmas \ref{ah} and \ref{lem:conseq-conn} respectively.
Since $g$ has infinite order, there exists $N_1\in \NN$ such that $g^n \notin \mathcal{F} \cup \Xi$ for any $n \in \Z$
with $|n| \ge N_1$.

Choose an arbitrary $n\in \Z$ such that $|n| \ge N_1$. By Lemma \ref{ah}, the element $g^nx=(xg^n)^x$ is loxodromic in $G$.
Suppose that $y \in E_G(g^nx)$. By part (2) of Lemma \ref{Eg}, there are $m \in \NN$ and $\epsilon \in \{-1,1\}$ such that
\begin{equation}\label{eq:m} y(g^nx)^my^{-1}=(g^nx)^{\epsilon m}.\end{equation}

Let $V$ be the letter from $\cH'$ representing $g^n$ in $G$, let $W$ be the letter from $X$ representing $x$, and let
$U$ be the shortest word over the alphabet $X \cup \cH'$ representing $y$.
Set $C=\|U\|$ and $d=1$. Now we apply Lemma \ref{lem:conseq-conn}.(b) to find the constant $L=L(C,d)$ from its claim.
Evidently we can assume that the number $m$ from equation \eqref{eq:m} is larger than $L$.

Consider a cycle $o=rqr'q'$ in $\Gamma(G,X\cup\cH')$ where
$\phi(r)\equiv U$, $\phi(q) \equiv (VW)^m$, $\phi(r')\equiv U^{-1}$,
$\phi(q') \equiv (VW)^{-\epsilon m}$. By construction, the cycle $o$
satisfies the assumptions of Lemma \ref{lem:conseq-conn}.(b), hence
some components $p$ of $q$ and $p'$ of $q'^{-1}$ must be connected
in $\Gamma(G,X\cup\cH')$. That is, there is a path $s$ with
$s_-=p_+$, $s_+=p'_+$ and $z=\phi(s)\in E_G(g)$ (see Figure \ref{pic:1}). Note that $\phi(p)
\equiv V$, $\phi(p')\equiv V^{\epsilon}$.

\begin{figure}[!ht]
  \begin{center}
   \input{pic1.pstex_t}
  \end{center}
  \caption{}\label{pic:1}
\end{figure}

Let $q_1$ be the subpath of $q$ starting at $r_+=q_-$ and ending at
$p_+=s_-$; let $q_1'$ be the subpath of $q'$ starting at $s_+=p'_+$
and ending at $q'_+=r_-$. Considering the cycle $o_1=rq_1sq_1'$ in
the case when $\epsilon=-1$ we get the following equality in $G$:
$$(g^nx)^\xi y (g^{n}x)^\zeta= z^{-1} g^{-n} \in E_G(g^nx) \cap E_G(g)~\mbox{ for some } \xi,\zeta \in \Z.$$
Similarly, in the case when $\epsilon=1$, we get
$$(g^nx)^\xi y (g^{n}x)^\zeta= g^n z^{-1} g^{-n} \in E_G(g^nx) \cap E_G(g)~\mbox{ for some } \xi,\zeta \in \Z.$$

Observe that by Lemma \ref{Eg}, the group $G$ is hyperbolic relatively to $\Hl \cup \{E_G(g),E_G(g^nx)\}$, hence, by
Lemma \ref{maln}, the intersection $E_G(g^nx) \cap E_G(g)$ is finite. Since $g$ is $H$-special, any finite subgroup of
$E_G(g)$ is contained in $E_G(H)$. Therefore $E_G(g^nx) \cap E_G(g) \subset E_G(H)$.
Thus, whatever $\epsilon \in \{-1,1\}$ is, we always have
$(g^nx)^\xi y (g^{n}x)^\zeta= h \in E_G(H)$, implying that $y=(g^nx)^{-\xi-\zeta}h$ because $g,x \in C_H(E_G(H))$.
By part (2) of Lemma \ref{Eg}, $\langle g^n x \rangle$ and $E_G(H)$ are both contained in $E_G(g^nx)$; consequently
$E_G(g^nx)=\langle g^n x \rangle \times E_G(H)$.
\end{proof}

\begin{proof}[Proof of Proposition \ref{fi}.] By Lemma \ref{SG} we can find an element $g \in S_G(H)$.
Note that for any $x \in Z= E_G(H) \cap C_H(E_G(H))$, the element
$gx$ is also $H$-special. Since $x=g^{-1}(gx)$, we have $Z \subset
\langle S_G(H) \rangle$. It is easy to see that $E_G(g)\cap
C_H(E_G(H)) = \langle g \rangle \times Z$, hence $E_G(g)\cap
C_H(E_G(H))\subset \langle S_G(H) \rangle$. Now, if $x \in
C_H(E_G(H))\setminus E_G(g)$, then by Lemma \ref{lem:spec-mod},
$g^nx \in S_G(H)$ for some $n \in \NN$. Consequently, $x=g^{-n} (g^n
x) \in \langle S_G(H) \rangle$.
\end{proof}


\section{Technical lemmas}


Our main goal here is to prove several auxiliary lemmas, which will be used in the next section to give
an algebraic description of automorphisms preserving commensurability classes of
elements in relatively hyperbolic groups. We begin with a definition.

\begin{defn}\label{def:commensurability}
Let $G$ be a group. Two elements $x,y \in G$, are said to be {\it commensurable} if there are $z \in G$, $m,n \in \Z \setminus \{0\}$
such that $y^n=zx^mz^{-1}$ in $G$. If the elements $x$ and $y$ are commensurable in $G$, we will write $x \stackrel{G}{\approx} y$; otherwise, we will write $x \napg y$.
\end{defn}

\begin{rem}\label{rem:lox-comm} Obviously any two elements of finite order are commensurable. Further, if $g$ and $h$ are commensurable elements of a relatively hyperbolic group $G$ and $g$ is loxodromic, then $h$ is loxodromic too. Indeed, evidently $h$ has infinite order. Suppose that $h$ is parabolic. Since $g \apg h$,
there are $\lambda \in \Lambda$, $a \in G$ and $m \in \NN$ such that $a^{-1} g^m a \in H_\lambda$. Since $g$ is loxodromic,
$x=g^a \notin H_\lambda$ and the intersection $H^x \cap H$ contains an infinite order element $x^m$. The
latter contradicts claim (2) of Lemma \ref{maln}.
\end{rem}

Throughout the rest of this section, $G$ will denote a group hyperbolic relative to a
collection of peripheral subgroups $\Hl$, and $H \le G$ will denote a non-elementary subgroup with $H^0 \neq \emptyset$.
\begin{lem}\label{lem:non-comm} Let $g \in G$ be a loxodromic element and $x\in G \setminus E_G(g)$.
For any finite subset $Y$ of  $G$
there is $N_2 \in \NN$ such that $g^n x$ is loxodromic and is not commensurable with any $y \in Y$ whenever $|n| \ge N_2$.
\end{lem}

\begin{proof} In view of Lemma \ref{Eg}.(3), we can assume that $E_G(g)$ belongs to the family of peripheral
subgroups $\Hl$ of $G$ and each infinite order element $y \in Y$ is parabolic.

Now we can apply Lemma \ref{ah}, to find $N_2 \in \NN$ such that for any $n \in \Z$ with $|n|\ge N_2$, the element $xg^n$
is loxodromic.
Therefore, so is $h=g^nx=x^{-1}(xg^n)x$. Suppose that $h$ is commensurable with some $y \in Y$. Then $y$ must also be loxodromic
(by Remark \ref{rem:lox-comm}), which contradicts our assumption above.
\end{proof}

\begin{lem}\label{lem:general} Let  $\{g_1,\dots,g_l\}$, $l\ge 2$, be a finite set of pairwise non-commensurable loxodromic
elements in a relatively hyperbolic group $G$. For any finite subset $F \subset G$ there exists $N_3 \in \NN$ such that
for any permutation $\sigma$ of $\{1,2,\dots,l\}$ and arbitrary elements $h_i \in E_G(g_{\sigma(i)})$, $i=1,2,\dots,l$,
of infinite order, the following hold.
\begin{itemize}
\item[(i)] The element $g=g_1^{m_1}f_1g_2^{m_2}f_2\dots g_l^{m_l}f_l$ is loxodromic for any $f_i \in F$ and  $m_i \in \Z$ with
$|m_i|\ge N_3$, $i=1,2,\dots,l$.
\item[(ii)] Suppose that  $\left( g_1^{m_1}g_2^{m_2}\dots g_l^{m_l} \right)^\zeta$ is conjugate to
$\left( h_1^{n_1}f_1 h_2^{n_2} f_2 \dots h_l^{n_l}f_l\right)^\eta$ in $G$,
for some  $f_i \in F$, $\zeta,\eta \in \NN$, $m_i,n_i\in \Z$, $|m_i|\ge N_3$, $|n_i|\ge N_3$ for all $i=1,2,\dots,l$.
Then $\zeta=\eta$, there is $k\in \{0,1,\dots,l-1\}$ such that
$\sigma$ is a cyclic shift by $k$, that is $\sigma(i) \equiv i+k ~ ({\rm mod}~ l)$ for all $i \in \{1,2,\dots,l\}$, and
$f_j \in E_G \left( g_{\sigma(j)} \right) E_G\left(g_{\sigma(j+1)} \right)$ when $j=1,2,\dots,l-1$,
$f_l \in E_G\left( g_{\sigma(l)}\right) E_G \left( g_{\sigma(1)}\right)$.
\end{itemize}
\end{lem}

\begin{proof} By Lemma \ref{Eg} and because $g_i \napg g_j$ when $i\neq j$,
$G$ is hyperbolic relative to the extended collection of subgroups
$\Hl \cup \{E_G(g_i)\}_{i=1}^l$.
Also, the finite relative generating set $X$ can be replaced by the bigger finite
set $X'=X \cup F\cup F^{-1}$ retaining the relative hyperbolicity of $G$.
Denote
$\mathcal{H}'=\left(\cup_{\lambda\in \Lambda} H_\lambda \cup \cup_{i=1}^l E_G(g_i)\right)\setminus \{1\} \subset G$.
Let $\Xi$ be the finite subset of $G$ given by Lemma \ref{lem:conseq-conn}.

Take any $i \in \{1,\dots,l\}$. By part (1) of Lemma \ref{Eg}, we have $|E_G(g_i): \langle g_i \rangle|<\infty$, hence any infinite order element $h \in E_G(g_i)$
belongs to the elementary subgroup $$E^+_G(g_i)=\{x \in G~|~\exists~m \in \mathbb{N}~\mbox{such that } x^{-1}g_i^mx=g_i^m\} \le E_G(g_i).$$
Clearly, the center of $E^+_G(G_i)$ has finite index in it, hence all finite order elements of $E_G^+(g_i)$
form the maximal torsion subgroup $T \lhd E_G^+(g_i)$. Let $\alpha: E_G^+(g_i) \to E_G^+(g_i)/T$ be the natural epimorphism. The image
$\alpha(E_G^+(g_i))$ is infinite cyclic (because it is virtually cyclic and torsion-free), therefore there exists
 $K_i \in \NN$ such that for every non-trivial element $y \in \alpha(E_G^+(g_i))$, one has $y^{n} \notin S_i$ whenever $|n|\ge K_i$,
where $S_i=\alpha(E_G^+(g_i)\cap \Xi)$ is a finite subset of $\alpha(E_G^+(g_i))$.  Set $N_3=\max\{K_i~|~i=1,\dots,l\}$.
By construction, for every $i$ and each infinite order element $h \in E_G(g_i)$, we have
$h^{n} \notin \Xi$ whenever $|n| \ge N_3$.

Choose any elements $f_i \in F$ and integers $m_i$ with $|m_i|\ge N_3$, $i=1,\dots,l$.
Let $V_i$ and $W_i$ be the letters from $\cH'$ and from $X'$ representing the elements
$g_i^{m_i}$ and $f_i$, $i=1,\dots,l$, respectively.

Proving claim (i) by contradiction, suppose that the element $g$ is not loxodromic.

If $g$ has finite order $t \in \NN$,
then set $C=0$, $d=1$ and choose $L=L(C,d)$ according to Lemma \ref{lem:conseq-conn}.(b).
In the Cayley graph $\Gamma(G,X'\cup\cH')$ consider the cycle $o=rqr'q'$, where
$\phi(q) \equiv(V_1W_1V_2W_2\dots V_lW_l)^{Lt}$, and $r,r'$ and $q'$ are trivial
paths consisting of single vertex $q_-=q_+=1$.
Since $\L(q) \ge Lt \ge L$, it follows from Lemma \ref{lem:conseq-conn}.(b) that some component of $q$ must be connected to a component
of $q'^{-1}$. But $q'^{-1}$ has no components at all. A contradiction.

Therefore $g$ must have infinite order and must be parabolic, i.e., $g=aha^{-1}$ for some $h \in \cH'$ and $a\in G$.
Let $C=|a|_{X'\cup \cH'}$, $d=2$ and $L=L(C,d)$ be given by Lemma \ref{lem:conseq-conn}.(b).
Since $h$ has infinite order (as a conjugate of $g$), there is $n \in \NN$ such that $n \ge L$ and $h^n \notin \Xi$.
Choose a shortest word $A$ over $X'\cup \cH'$ representing $a$ in $G$, and let $U$ be the letter from
$\cH'$ corresponding to $h^n$. Consider a cycle $o=rqr'q'$
in $\Gamma(G,X'\cup\cH')$ such that $\phi(r)\equiv A$, $q_-=r_+$, $\phi(q)\equiv (V_1W_1V_2W_2\dots V_lW_l)^{n}$,
$r'_-=q_+$, $\phi(r')\equiv A^{-1}$, $q'_-=r'_+$, $\phi(q')\equiv U^{-1}$. Since $\L(r)=\L(r')=C$, $\L(q)\ge n\ge L$, we can
apply Lemma~\ref{lem:conseq-conn}.(b) to $o$, claiming that two distinct components of $q$ must be connected to two
distinct components of $q'^{-1}$. But $q'^{-1}$ has only one component by definition.
This contradiction concludes the proof of claim (i).

To establish claim (ii), assume that
$ b\left( g_1^{m_1}g_2^{m_2}\dots g_l^{m_l} \right)^\zeta b^{-1}=
\left( h_1^{n_1}f_1 h_2^{n_2} f_2 \dots h_l^{n_l}f_l\right)^\eta$ in $G$,
for some infinite order elements $h_i \in E_G(g_{\sigma(i)})$,  $b \in G$, $\zeta,\eta \in \NN$, $m_i,n_i\in \Z$, $|m_i|\ge N_3$, $|n_i|\ge N_3$ for $i=1,2,\dots,l$.
Then for every $\varkappa \in \NN$ we have
\begin{equation}\label{eq:b-zeta} b\left( g_1^{m_1}g_2^{m_2}\dots g_l^{m_l} \right)^{\varkappa \zeta} b^{-1}=
\left( h_1^{n_1}f_1 h_2^{n_2} f_2 \dots h_l^{n_l}f_l\right)^{\varkappa \eta} .
\end{equation}
Let $V_i$ and $W_i$ be as before. Choose a letter $U_i$ from $\cH'$ corresponding to $h_{i}^{n_i}$, $i=1,\dots,l$,
and a shortest word $B$ over $X'\cup \cH'$ representing $b$ in $G$.
Set $C=\|B\|$, $d=2l$ and let $L=L(C,d)\in \NN$ be the constant given by Lemma \ref{lem:conseq-conn}.(b).
Take $\varkappa \in \NN$ so that $\varkappa\zeta l \ge L$ and $\varkappa l >6C$.

In the Cayley graph $\Gamma(G,X'\cup\cH')$ equation \eqref{eq:b-zeta} gives rise to a cycle $o=rqr'q'$, in which
$\phi(r)\equiv B$, $q_-=r_+$, $\phi(q)\equiv (V_1V_2\dots V_l)^{\kappa\zeta}$,
$r'_-=q_+$, $\phi(r')\equiv B^{-1}$, $q'_-=r'_+$,
$\phi(q')\equiv \left( U_1 W_1 U_2 W_2\dots U_l W_l \right)^{-\kappa\eta}$.

By construction, the paths $q$ and $q'$ have exactly $\varkappa \zeta l$ and $\varkappa\eta l$ components respectively.
Suppose that $\zeta >\eta$. By Lemma \ref{lem:conseq-conn}.(a), at least
$\varkappa\zeta l-6C > \varkappa l(\zeta-1)\ge \varkappa l\eta$ components of
$q$ must be connected to components of $q'$, hence two distinct components of $q$ will have to be connected
to the same component of $q'$, contradicting Lemma \ref{lem:conseq-conn}.(a). Hence $\zeta \le \eta$. A symmetric argument
shows that $\eta \le \zeta$. Consequently $\zeta=\eta$.

Since $\L(q)=\varkappa \zeta l \ge \L$, we can apply Lemma
\ref{lem:conseq-conn}.(b) to find $2l$ consecutive components of $q$
that are connected to $2l$ consecutive components of $q'^{-1}$.
Therefore there are consecutive components $p_1,\dots ,p_{l+1}$ of
$q$ and $p'_1,\dots, p'_{l+1}$ of $q'^{-1}$ such that $p_j$ is
connected to $p'_j$ for each $j$, and $\phi(p_i)\equiv V_i$ for
$i=1,\dots,l$, $\phi(p_{l+1})\equiv V_1$
(Figure \ref{pic:2}). Therefore $\phi(p_i') \in
E_G(g_i)$, $i=1,\dots,l$, $\phi(p_{l+1}') \in E_G(g_1)$. From the
form of $\phi(q'^{-1})$ it follows that there is $k \in
\{0,1,\dots,l-1\}$ such that $\phi(p_j')\equiv U_{j+k}$ for
$j=1,\dots,l+1$  (indices at $U$ are taken modulo $l$). Thus
$U_{j+k}= h_{j+k}^{n_{j+k}}\in E_G(g_j)$. On the other hand,
$h_{j+k}^{n_{j+k}} \in E_G(g_{\sigma(j+k)})$ has infinite order.
Hence the intersection $E_G(g_j) \cap E_G(g_{\sigma(j+k)})$ must be
infinite, which yields (by Lemma \ref{maln}) that $\sigma(j+k)=j$ for
all $j$. Therefore $\sigma$ is a cyclic shift (by $l-k$) of
$\{1,\dots,l\}$.

\begin{figure}[!ht]
  \begin{center}
   \input{pic2.pstex_t}
  \end{center}
  \caption{}\label{pic:2}
\end{figure}

The subpath $w_i$ of $q'^{-1}$ between
$(p'_i)_+$ and $(p'_{i+1})_-$ is labelled by $W_{\sigma^{-1}(i)}$. As we showed, the vertex $(p_i)_+=(p_{i+1})_-$
is connected to $(w_i)_-$ by a path $s_i$ with $\phi(s_i) \in E_G(g_i)$, and  to $(w_i)_+$
by a path $t_i$ with $\phi(t_i) \in E_G(g_{i+1})$, $i =1,\dots,l$ (here we use the convention that $g_{l+1}=g_1$). Considering the cycle
$t_i^{-1}s_i w_i$ we achieve the desired inclusion: $f_{\sigma^{-1}(i)}=\phi(w_i) \in E_G(g_i) E_G(g_{i+1})$,
$i=1,\dots,l$. This concludes the proof.
\end{proof}

\begin{lem}\label{lem:prod_of_three} Suppose that $\varphi: H \to G$ is a homomorphism such that
$\varphi(h) \apg h$ for all $h \in H^0$. Then for any $g_1,g_2,g_3 \in H^0$,
satisfying $g_i \stackrel{G}{\not\approx} g_j$
for $i\neq j$,
there exists $N_4 \in \NN$ such that for arbitrary $n_1,n_2,n_3 \in \Z$, with $|n_i| \ge N_4$, $i=1,2,3$, and for
$g =g_1^{n_1} g_2^{n_2} g_3^{n_3}$, one has $g \in H^0$ and
$(\varphi(g))^\zeta= e g^\zeta e^{-1}$, for some $e \in G$ and $\zeta \in \NN$.
\end{lem}

\begin{proof} According to the assumptions, there exist $x_i \in G$ and $\zeta_i,\eta_i \in \Z\setminus\{0\}$ such that
$\left( \varphi(g_i) \right)^{\zeta_i}= x_i g_i^{\eta_i} x_i^{-1}$, $i=1,2,3$. Denote $h_i=x_i^{-1} \varphi(g_i) x_i$, $i=1,2,3$.
Then $h^{\zeta_i}=g^{\eta_i}$ , hence $h_i \in E_G(g_i)$ (by part (2) of Lemma \ref{Eg}) and $h_i$ has infinite order, $i=1,2,3$.

Set $f_1=x_1^{-1}x_2$,
$f_2=x_2^{-1}x_3$ and $f_3=x_3^{-1}x_1$, and let $N_4 \in \NN$ be the number $N_3$ from the claim of
Lemma \ref{lem:general} applied to
the set of loxodromic elements $\{g_1,g_2,g_3\}$ and the finite set $F=\{f_1,f_2,f_3\}$.
Take any $n_i \in \Z$ with $|n_i| \ge N_4$, $i=1,2,3$.
By part (i) of Lemma \ref{lem:general}, $g =g_1^{n_1} g_2^{n_2} g_3^{n_3} \in H^0$. Hence there are
$\zeta,\eta \in \Z\setminus \{0\}$ and $e \in G$ such that $e g^\zeta e^{-1}=(\varphi(g))^\eta $.
Since $\varphi$ is a homomorphism, we get
$$ e(g_1^{n_1} g_2^{n_2} g_3^{n_3})^\zeta e^{-1}= (\varphi(g))^\eta =
(x_1h_1^{n_1}x_1^{-1} x_2 h_2^{n_2} x_2^{-1} x_3 h_3^{n_3} x_3^{-1})^\eta , ~\mbox{ hence }$$
\begin{equation}\label{eq:three} (x_1^{-1}e)(g_1^{n_1} g_2^{n_2} g_3^{n_3})^\zeta (x_1^{-1}e)^{-1}=
(h_1^{n_1}f_1 h_2^{n_2} f_2h_3^{n_3} f_3)^\eta .
\end{equation}

Without loss of generality we can assume that $\zeta >0$. Suppose that $\eta<0$. Then
$(g_3^{-n_3} g_2^{-n_2} g_1^{-n_1})^\zeta$ is conjugate to $(h_1^{n_1}f_1 h_2^{n_2} f_2h_3^{n_3} f_3)^{-\eta} $
in $G$ and $-\eta>0$.
Applying part (ii) of Lemma \ref{lem:general} to this situation, we get a contradiction with
the fact that the transposition $(1,3)$ is not
a cyclic shift of $\{1,2,3\}$. Therefore, $\eta>0$ and we can apply part (ii) of Lemma \ref{lem:general}
to \eqref{eq:three}, achieving the required equality $\zeta=\eta$.
\end{proof}

\begin{lem}\label{lem:first_step} Let $a,b \in G$ be non-commensurable loxodromic elements and let $y,z \in G$.
There exists $N_5\in \NN$ such that the following holds. Suppose that
$a^{k'} y b^{l'} z\stackrel{G}{\approx} a^{k} b^{l}$ for some integers $k,l,k',l'$ with
$|k|,|l|,|k'|,|l'|\ge N_5$. Then  $y \in E_G(a) E_G(b)$ and $z \in E_G(b) E_G(a)$.
\end{lem}

\begin{proof} Choose $N_5 \in \NN$ to be the number $N_3$ arising after an application of
Lemma \ref{lem:general} to $\{a,b\}$ and $F=\{y,z\}$. Choose any $k,l,k',l' \in \Z$ satisfying
$|k|,|l|,|k'|,|l'|\ge N_5$.

Assume that there is $e \in G$, $\zeta \in \NN$ and $\eta \in \Z\setminus \{0\}$
for which $e\left( a^{k} b^{l}\right)^\zeta e^{-1} = \left(a^{k'} y b^{l'} z\right)^\eta$. If $\eta>0$ then the statement
immediately follows from part (ii) of Lemma \ref{lem:general}. So, suppose that $\eta <0$. Then $-\eta >0$ and
$\left(  b^{-l}a^{-k}\right)^\zeta$ is conjugate to $\left(a^{k'} y b^{l'} z\right)^{-\eta}$ in $G$. Again, by part
(ii) of Lemma \ref{lem:general}, $y \in E_G(a) E_G(b)$, $z \in E_G(b) E_G(a)$.
\end{proof}

\begin{lem} \label{lem:spec-image} Assume that $g \in S_G(H)$ and $\psi: H \to G$ is a homomorphism satisfying
$\psi(g^n)=g^n z$ for some $n \in \NN$ and $z \in E_G(H)$. Then there is $f \in E_G(H)$ such that $\psi(g)=gf$.
\end{lem}

\begin{proof} After replacing $n$ with $n'=n|E_G(H)|$, we can further assume that $z=1$, because
$\psi(g^{n'})=g^{n'} z^{n'}=g^{n'}$.

Now, note that $\psi(g) g^n (\psi(g))^{-1}= \psi(g^n) =g^{n}$, hence $\psi(g) \in E_G(g)$
by part (2) of Lemma \ref{Eg}. Since $g$ is $H$-special, there is $k \in \Z$ and $f \in E_G(H)$ such that $\psi(g)=g^k f$.
Denote $l=|E_G(H)|$. Then $g^{ln}=\psi(g^{ln})=(g^k f)^{ln}=g^{lnk}f^{ln}=g^{lnk}$, implying that $k=1$, as required.
\end{proof}

\begin{lem}\label{lem:lox_non-comm} Suppose that for an automorphism $\alpha \in Aut(H)$ there is $g \in H^0$ satisfying $g \napg \alpha(g)$. Then there
exists an element $a \in H$ such that both $a$ and  $\alpha(a)$ are loxodromic in $G$ and $a \napg \alpha(a)$.
\end{lem}

\begin{proof} If $\alpha(g) \in H^0$, there is nothing to prove. Thus, we can assume that $\alpha(g)$ is parabolic in $G$, i.e., there exists a
peripheral subgroup $H_\lambda$ and elements $t\in G$, $h \in H_\lambda$ such that $\alpha(g) =h^t$.
Denote $x=\alpha^{-1}(g) \in H$. If $x \in E_G(g)$, then $\langle g \rangle^x \cap \langle g \rangle$ is infinite (by Lemma \ref{Eg}.(b)), hence
$\langle \alpha(g) \rangle^{\alpha(x)} \cap \langle \alpha(g) \rangle$ is infinite. Thus
$H_\lambda^{(tgt^{-1})} \cap H_\lambda$ is infinite, which implies, by Lemma \ref{maln}, that $tgt^{-1} \in H_\lambda$, contradicting the loxodromicity of $g$.

Therefore $x \notin E_G(g)$. Since both $g$ and $\alpha(g)$ have infinite order and $y=tgt^{-1} \in G\setminus H_\lambda$, we can apply Lemmas \ref{lem:non-comm}
and \ref{ah} to find $N \in \NN$ such that for any integer $n \ge N$, the elements $g^n x$ and $h^n y$ are loxodromic in $G$. Note that
$\alpha(g^n x)=(h^n y)^t$.

Suppose, first, that
\begin{equation} \label{eq:g^nx} g^n x \apg \alpha(g^n x)~\mbox{ for every } n \ge N.\end{equation}
 By Lemma \ref{Eg}, $G$ is hyperbolic relative to $\Hl \cup \{E_G(g)\}$.
Without loss of generality, we can also assume that $x$ and $y$ belong to the finite relative generating set $X$ of $G$. Let $\Xi \subset G$ be the finite
set from Lemma \ref{lem:conseq-conn}. Evidently there is an integer $n \ge N$ such that $g^n,h^n \notin \Xi$. Our assumption \eqref{eq:g^nx} implies that there is $b \in G$,
$k,l \in \Z\setminus \{0\}$ such that $b (g^n x)^k b^{-1}=(h^n y)^l$. Choose a word $B$ in the alphabet
$X \cup \mathcal{H}'$ representing $b$ in $G$, where $\mathcal{H}'=(\cup_{\lambda\in \Lambda} H_\lambda  \cup E_G(g))\setminus\{1\}$, and
let $W,Y \in X$, $U \in E_G(g)$, $V \in H_\lambda$ be the letters corresponding to $x,y,g^n,h^n$ respectively. Set $d=1$, $C=\|B\|$ and let $L=L(C,d)$ be the constant
provided by part (b) of Lemma \ref{lem:conseq-conn}. Without loss of generality we can assume that $|k|,|l| \ge L$.

Consider a cycle $o=rqr'q'$ in the Cayley graph  $\Gamma(G,X\cup\mathcal{H}')$, where $\phi(r) \equiv B$, $r_+=q_-$, $\phi(q) \equiv (UW)^k$, $q_+=r'_-$,
$\phi(r') \equiv B^{-1}$, $q_-'=r'_+$ and
$\phi(q') \equiv (VY)^{-l}$. It is easy to see that $o$ satisfies all the conditions of Lemma \ref{lem:conseq-conn}, hence
some component of $q$ must be connected to a component of $q'^{-1}$ in $\Gamma(G,X\cup\mathcal{H}')$. However,
according to the construction, $q$ has only $E_G(g)$-components, and $q'^{-1}$ has only $H_\lambda$-components. Thus the assumption
\eqref{eq:g^nx} yields a contradiction. Hence, there exists $n \ge N$ such that for the element $a=g^nx$ 
we have $a \in H^0$, $\alpha(a) \in H^0$ and $a \napg \alpha(a)$.
\end{proof}


\section{Commensurating automorphisms of relatively hyperbolic groups}\label{sec:comm-aut}


The purpose of this section is to study automorphisms of relatively hyperbolic groups preserving commensurability classes. 

Recall that $N_G(H)$ denotes the normalizer of a subgroup $H$ in a group $G$. Further, let $H$ be a non-elementary subgroup of a relatively hyperbolic group $G$ such that $H^0\ne \emptyset$. We denote by $\widehat H$ the product $H E_G(H)$. This is clearly a subgroup of $G$.

\begin{thm} \label{thm:comm-aut}
Let $G$ be a relatively hyperbolic group, let $H \le G$ be a
non-elementary subgroup and let $\varphi \in Aut(H)$.
Suppose that $H^0\ne \emptyset $ and $\varphi (h)\apg h$ for every $h\in H^0$.
Then there is a set map $\e : H \to E_G(H)$, whose restriction to $\CH$ is a homomorphism,
and an element $w \in N_G (\widehat H )$ such that for every $h \in H$,
$\varphi(h)=w \left( h \e (h) \right)w^{-1}$.
\end{thm}

Below is them main technical lemma of this section. It demonstrates
how to construct the element $w$ and the restriction of the map $\e$ to $\CH$ from the statement of Theorem \ref{thm:comm-aut}.

\begin{lem} \label{lem:comm-aut-centralizer} Suppose that $G$ is a relatively hyperbolic group, $H \le G$ is a
non-elementary subgroup and $\varphi \in Aut(H)$. Assume that $H^0\ne \emptyset $ and $\varphi (h)\apg h$ for every $h\in H^0$. Then there
is a homomorphism $\tilde \e : C_H \left( E_G(H) \right) \to E_G(H)$ and an element $w \in G$ such that for every $x \in C_H(E_G(H))$,
$\varphi(x)=w \left( x \tilde \e (x) \right)w^{-1}$.
\end{lem}

\begin{proof}
By Lemma \ref{SG}, $H$ contains an $H$-special element $g_1$. Since $H$ is non-elementary and $\CH$ has finite index in it, $\CH$
is also non-elementary. The subgroup $E_G(g_1)$ is elementary (by part (1) of Lemma \ref{Eg}), thus
there is an element $y \in \CH \setminus E_G(g_1)$.
By Lemma \ref{lem:non-comm}, there is $k_2 \in \NN$ such that $g_2=g_1^{k_2}y \in \CH$ is loxodromic and $g_2 \napg g_1$.
Using the same lemma again we can find $k_3 \in \NN$ such that $g_3=g_1^{k_3}y \in \CH$ is loxodromic and $g_3 \napg g_i$, $i=1,2$.
In particular, $E_G(g_2) \cap \langle g_3 \rangle = \{1\}$.

Choose $N_4 \in \NN$ according to an application of Lemma \ref{lem:prod_of_three} to $\varphi$, $g_1,g_2,g_3$, and let $n_3=N_4$.
By Lemma \ref{lem:non-comm}, there is $n_2 \ge N_4$ such that $g_2^{n_2}g_3^{n_3} \in H^0$ is not commensurable with $g_1$ in $G$.
Therefore $g_2^{n_2}g_3^{n_3}\in \CH \setminus E_G(g_1)$,
and by Lemma \ref{lem:spec-mod} there is $N_1 \in \NN$ such that the element
$g_1^{n}g_2^{n_2}g_3^{n_3}$ is $H$-special for any $n\ge N_1$.
Denote $n_1=\max\{N_1,N_4\}$ and apply Lemma \ref{lem:non-comm} to find $m \in \NN$ such that the elements
$a=g_1^{n_1}g_2^{n_2}g_3^{n_3}$ and
$b=g_1^{n_1+m}g_2^{n_2}g_3^{n_3}$ are not commensurable with each other in $G$.
In view of Lemma \ref{lem:prod_of_three} one can conclude that the elements $a,b \in \CH$ are
$H$-special and there exist $u,v \in G$, $\mu,\nu \in \NN$ such that $\varphi(a^\mu)=u a^\mu u^{-1}$,
$\varphi(b^\nu)=vb^\nu v^{-1}$.

Let $\chi: H \to G$ be the monomorphism, defined by $\chi(h)=u^{-1} \varphi(h) u$
for all $h \in H$. Then $\chi(a^\mu)=a^\mu$, $\chi(b^\nu)=(u^{-1}v)b^{\nu}(u^{-1}v)^{-1}$.
Note that $\chi(h) \apg h$ for every $h \in H^0$. By part (i) of Lemma \ref{lem:general},
$a^{k\mu}b^{k\nu} \in H^0$ for every sufficiently large $k \in \NN$. Therefore
 \begin{equation*}\label{eq:ab}
a^{k\mu} (u^{-1}v) b^{k\nu}(u^{-1}v)^{-1}=\chi(a^{k\mu}b^{k\nu}) \apg a^{k\mu} b^{k\nu}~\mbox{ for every sufficiently large } k \in \NN.
\end{equation*}
Consequently, by Lemma \ref{lem:first_step},
$u^{-1}v \in E_G(a)E_G(b)$, thus $u^{-1}v= a^s b^t f$ for some $s,t \in \Z$, $f \in E_G(H)$. Hence
$\chi(b^\nu)=a^sb^{\nu} a^{-s}$ because $b \in \CH$. Denote $w=ua^s \in G$ and let $\psi: H \to G$ be the monomorphism
defined by the formula $\psi(h)=w^{-1} \varphi(h)w=a^{-s} \chi(h) a^s$ for all $h\in H$.
By construction, we have
\begin{equation}\label{eq:psi} \psi(a^\mu)=a^\mu,~ \psi(b^\nu)=b^\nu ~\mbox{ and }~\psi(h) \apg h \mbox{ for each }h \in H^0.
\end{equation}

Choose any element $g \in S_G(H)$. We will show that there is $f \in E_G(H)$
such that  $\psi(g)=gf$.

If $g \in E_G(a)$ then there is $n \in \NN$ such that
$g^n \in \langle a^\mu \rangle$ because $|E_G(a):\langle a^\mu \rangle |<\infty$. Hence $\psi(g^n)=g^n$ and by
Lemma \ref{lem:spec-image}, $\psi(g)=gf$ for some $f \in E_G(H)$.

Suppose, now, that $g \notin E_G(a)$. Since $g \in \CH$ and $a$ is $H$-special, we can use
Lemmas \ref{lem:spec-mod} and \ref{lem:non-comm} to find $l \in \NN$ such that the element $d=a^{l\mu}g$ is $H$-special and is not
commensurable with $a$ and $b$ in $G$. Arguing as in the beginning of the proof
(using Lemmas \ref{lem:spec-mod}, \ref{lem:non-comm} and \ref{lem:prod_of_three})
we can find $m_1,m_2,m_3 \in \NN$ such that
$c=a^{m_1 \mu} b^{m_2 \nu} d^{m_3} \in S_G(H)$,  $c \napg a$, $c \napg b$ and $\psi(c^\zeta)=e c^\zeta e^{-1}$ for some $\zeta \in \NN$
and $e \in G$.

By part (i) of Lemma \ref{lem:general}, $a^{k\mu}c^{k\zeta} \in H^0$ for every sufficiently
large $k \in \NN$. Hence $a^{k\mu} e c^{k \zeta} e^{-1}=\psi\left( a^{k\mu} c^{k\zeta} \right) \apg a^{k\mu} c^{k\zeta}$ whenever
$k$ is sufficiently large. Applying Lemma \ref{lem:first_step} we see that $e \in E_G(a)E_G(c)$. As before, this implies that
$\psi(c^\zeta)= a^p c^\zeta a^{-p}$ for some $p \in \Z$.

Similarly, there is $q \in \Z$ such that $\psi(c^\zeta)= b^q c^\zeta b^{-q}$. Hence
$(a^{-p}b^q) c^\zeta (a^{-p}b^q)^{-1}=c^\zeta$, yielding that $a^{-p}b^q \in E_G(c)$.

Suppose that $p \neq 0$ and $q\neq 0$. Then the element $a^{-p}b^q$ must have infinite order (otherwise we would have
$a^{-p}b^q \in E_G(H)$ since $c$ is $H$-special, hence $b^q \in a^p E_G(H) \subset E_G(a)$ contradicting to $a \napg b$).
This implies that $(a^{-p}b^q)^\alpha=c^\beta$ for some $\alpha\in \Z\setminus \{0\}$ and $\beta \in \NN$.
Recalling \eqref{eq:psi}, we can apply Lemma \ref{lem:spec-image} to find $f_1,f_2 \in E_G(H)$ such that $\psi(a)=af_1$
and $\psi(b)=bf_2$. Since $a,b \in \CH$ we obtain
$$\psi(c^\beta)=\psi \left((a^{-p}b^q)^\alpha\right)=\left(a^{-p}b^{q}\right)^\alpha f_3=c^\beta f_3~
\mbox{ for some } f_3 \in E_G(H).$$
Then for $\gamma=\beta\zeta |E_G(H)|$ we get $c^\gamma=\psi(c^\gamma)=a^p c^\gamma a^{-p}$,
implying that $a^p \in E_G(c)$, which contradicts to $a \napg c$.

Therefore either $p=0$ or $q=0$, thus $\psi(c^\zeta)=c^\zeta$. By Lemma \ref{lem:spec-image}, there is $f_5 \in E_G(H)$ such
that $\psi(c)=cf_5$. Since $c=a^{m_1 \mu} b^{m_2 \nu} d^{m_3}$, we can use \eqref{eq:psi} to get
$\psi(d^{m_3})=d^{m_3}f_5$. Applying Lemma
\ref{lem:spec-image} again, we find $f_6 \in E_G(H)$ such that $\psi(d)=df_6$. Finally, since $d=a^{l\mu}g$,
in view of \eqref{eq:psi} we achieve $\psi(g)=g f_6$, as needed.

To finish the proof, we observe that by Proposition \ref{fi}, $\CH$ is generated by $S_G(H)$, therefore for each
$x \in \CH$ there is $\tilde \e(x) \in E_G(H)$ such that $\psi(x)=x \tilde \e(x)$. Since $\psi$ is a homomorphism, the map
$\tilde \e:\CH \to E_G(H)$ will be a homomorphism too. By construction, we have
$\varphi(x)=w \psi(x) w^{-1}= w x \tilde \e(x) w^{-1}$.
\end{proof}

Now we are ready to prove the main result of this section.
\begin{proof}[Proof of Theorem \ref{thm:comm-aut}] Let $w \in G$ and $\tilde \e:\CH \to E_G(H)$ be as in the
claim of Lemma \ref{lem:comm-aut-centralizer}. Let $\psi:H \to G$ be the monomorphism that is defined according to the formula
$\psi(h)=w^{-1} \varphi(h) w$ for all $h \in H$. Denote $l=|H:\CH|$,  $m=|E_G(H)|$ and $n=ml \in \NN$.

Since $\CH$ is a normal subgroup of $H$, for any $z \in H$ we
have $z^l \in \CH$ and $\psi(z^n)=z^n \tilde\e(z^l)^m=z^n$. Fix an arbitrary $h \in H$. For any $g \in H^0$ we see that $g^n,hg^nh^{-1} \in \CH \cap H^0$, therefore
$\psi(h) g^n \psi(h)^{-1}=\psi(hg^nh^{-1})=hg^nh^{-1}$, implying that $h^{-1}\psi(h) \in E_G(g)$. Thus,
$h^{-1}\psi(h) \in \bigcap_{g\in H^0} E_G(g)=E_G(H)$. After defining $\e(h)=h^{-1}\psi(h)$ for each $h \in H$, one immediately sees that $\e:H \to E_G(H)$ is a
map with the required properties. Obviously, the restriction of $\e$ to $\CH$ coincides with $\tilde \e$.

It remains to prove that $w \in N_G(\widehat H)$. We will first show that $w \in N_G(E_G(H))$.
Consider any element $f \in E_G(H)$. Since  $\varphi$ is an automorphism of $H$, for any $g \in H^0$ there is $h \in H$ such that $\varphi(h)=g$.
Then $h^n \in \CH$ and $g^n=\varphi(h^n)=wh^nw^{-1}$ because $\e(h^n)=\tilde\e(h^l)^m=1$. Now we observe that
$$wfw^{-1} g^n (wfw^{-1})^{-1}=w f h^n f^{-1} w^{-1}=wh^nw^{-1}=g^n.$$
Hence, $wfw^{-1} \in E_G(g)$ for every $g \in H^0$; consequently $w f w^{-1} \in E_G(H)$. The latter implies that $w E_G(H) w^{-1} \subseteq E_G(H)$ and since
$E_G(H)$ is finite, we conclude that $w \in N_G(E_G(H))$.

Now, for any $h \in H$ we have $$whw^{-1}=w h \e(h)w^{-1} w \e(h)^{-1} w^{-1} = \varphi(h) \left( w \e(h) w^{-1} \right)^{-1}  \in H E_G(H);$$
thus $w H w^{-1} \subseteq \widehat H$. Since $w^{-1} \varphi(h) w = h \e(h) \in H E_G(H)$ and $\varphi \in Aut(H)$, one gets $w^{-1} H w \subseteq \widehat H$.
Therefore $w \widehat H  w^{-1} \subseteq \widehat H w E_G(H) w^{-1}=\widehat H$, $w^{-1} \widehat H  w \subseteq \widehat H w^{-1} E_G(H) w=\widehat H$, i.e., $w \in N_G(\widehat H)$.
\end{proof}

We are now in a position to prove Corollary \ref{cor:descr_comm_aut} mentioned in the Introduction.
We establish it in a more general form:

\begin{cor}  \label{cor:comm_aut_def}
Let $G$ be a non-elementary relatively hyperbolic group and $\varphi \in Aut(G)$. The following conditions are equivalent:
\begin{enumerate}
\item[(a)] $\varphi $ is commensurating;
\item[(b)] $\varphi (g)\apg g$ for every loxodromic $g\in G$;
\item[(c)] there is a set map $\e : G \to E(G)$, whose restriction to $C(G)$ is a homomorphism,
and an element $w \in G$ such that for every $g \in G$, $\varphi(g)=w \left( g \e (g) \right)w^{-1}$.
\end{enumerate}

In particular, if $E(G)=\{ 1\}$, then every commensurating automorphism of $G$ is inner.
\end{cor}

\begin{proof} (a) implies (b) by definition, and (b) implies (c) by Theorem \ref{thm:comm-aut}.
It remains to show that (c) implies (a). Indeed, let $g$ be an arbitrary element of $G$, and let the automorphism $\varphi $ satisfy (c).
If $g$ is of finite order, then so is $\varphi(g)$, and in this case evidently $\varphi (g)\apg g$.
Thus, we can suppose that $g$ has infinite order in $G$. By our assumptions,
$\varphi (g)=w(g\e(g))w^{-1}$ for some $w\in G$ and $\e(g)\in E(G)$. Since $E(G)$ is finite and normal in $G$, $\langle g\rangle $ has finite index in the subgroup $\langle g\rangle E(G)$. Hence there exists a non-zero integer $k$ such that $(g\e (g))^k=g^l$ for some $l \in \Z$. And since
the order of $g \e(g)=w^{-1} \varphi(g) w$ is infinite, we can conclude that $l\neq 0$.
Therefore $\varphi(g)=wg\e (g)w^{-1}$ is commensurable with $g$ in $G$. Thus $\varphi $ in commensurating.
\end{proof}



Recall that a result of Metaftsis and Sykiotis \cite[Lemma $2.2'$]{MS} states that for any relatively hyperbolic group $G$, one
has $|Aut_c(G):Inn(G)|<\infty$, where
$$Aut_c(G)=\{\alpha \in Aut(G)~|~\forall \, g\in G~\exists\, x=x(g) \in G~ \mbox{ such that } \alpha(g)=xgx^{-1}\}$$ is the group
of all {\it pointwise inner automorphisms} of $G$. Theorem \ref{thm:comm-aut} allows one to generalize their result to
all non-elementary subgroups:

\begin{cor}\label{cor:conj_aut} Suppose that $H$ is a non-elementary subgroup of a relatively
hyperbolic group $G$, with $H^0 \neq\emptyset$. Then $|Aut_c(H):Inn(H)|<\infty$. If, in addition, $E_G(H)=\{1\}$, then
$Aut_c(H)=Inn(H)$.
\end{cor}

\begin{proof} By Theorem \ref{thm:comm-aut}, for any automorphism $\varphi \in Aut_c(H)$,
there exist $w \in G$ and a map $\e:H \to E_G(H)$ such that $\varphi(h)=wh \e(h) w^{-1}$ for each $h \in H$.
Take any element $h \in S_G(H)$. Then $h$ commutes with $\e(h) \in E_G(H)$, and, consequently,
$(\varphi(h))^n=wh^nw^{-1}$ where $n=|E_G(H)|\in \NN$.

Now, since $\varphi$ is a pointwise inner automorphism of $H$, there
is $x \in H$ such that $\varphi(h)=xhx^{-1}$. Hence $x h^n x^{-1}=wh^nw^{-1}$, i.e.,
$w^{-1}x \in E_G(h)=\langle h \rangle \times E_G(H)$. Thus
$w = f z$ for some $f \in H$ and $z \in E_G(H)$, and $w^{-1}x \in C_G(h)$
because $h$ is $H$-special. Consequently, we have $h=w^{-1}x h \left(w^{-1}x\right)^{-1}=h \e(h)$, which implies that $\e(h)=1$.
Since the latter holds for any $h \in S_G(H)$, it follows from Proposition \ref{fi} that $\e(C_H)=\{1\}$, where $C_H=\CH$.

Note that $|H:C_H|<\infty$, hence there are $h_1,\dots,h_l \in H$ such that $H=\bigsqcup_{i=1}^l C_H h_i$.
For any $g \in H$ there are $a \in C_H$ and $i \in \{1,\dots, l\}$ such that $g=a h_i$. One has
\begin{multline*}\varphi(a) \varphi(h_i) = \varphi(g)=w g \e(g) w^{-1}=waw^{-1} w h_i \e(a h_i) w^{-1}= \\ \varphi(a)
\varphi(h_i) w (\e(h_i))^{-1}\e(a h_i) w^{-1}, \end{multline*}
hence $\e(g)=\e(a h_i)=\e(h_i)$, i.e.,
the map $\e$ is uniquely determined by the images of $h_1,\dots,h_l$.
Thus, $\varphi(g)=f z (g \e(h_i)) z^{-1} f^{-1}$, implying that the automorphism
$\varphi \in Aut_c(H)$, up to composition with an inner automorphism of $H$,
is completely determined by the finite collection of elements $z,\e(h_1),\dots,\e(h_l) \in E_G(H)$, and since
$E_G(H)$ is finite, we can conclude that $|Aut_c(H):Inn(H)|<\infty$.

Now, if $E_G(H)=\{1\}$ we obtain $w =f \in H$ and $\varphi(g)=wgw^{-1}$ for all $g \in H$, that is $\varphi \in Inn(H)$.
\end{proof}


\section{Group-theoretic Dehn surgery and normal automorphisms}\label{sec:Dehn_surgery}


In the context of relatively hyperbolic groups, the algebraic
analogue of Dehn filling is defined as follows. Suppose that $\Hl $
is a collection of (peripheral) subgroups of a group $G$. To each collection $\N
=\{ N_\lambda \} _{\lambda \in \Lambda }$, where $N_\lambda $ is a
normal subgroup of $H_\lambda $, we associate the quotient-group
\begin{equation}
G(\N ) = G/\ll \mbox{$\bigcup_{\lambda \in \Lambda }$} N_\lambda
\rr^G .
\end{equation}

\begin{defn}\label{def:periph_fill}
Let $G$ and $\Hl$ be as described above. We say that some assertion
holds for {\it most peripheral fillings of $G$}, if there exists a
finite subset $\mathcal F$ of non-trivial elements of $G$ such that
the assertion holds for $G(\N )$ for any collection $\N =\{
N_\lambda \} _{\lambda \in \Lambda }$ of normal subgroups $N_\lambda
\lhd H_\lambda $ satisfying $N_\lambda \cap \mathcal F=\emptyset $
for all $\lambda \in \Lambda $.
\end{defn}

The theorem below was proved in \cite{CEP}. In the particular case
when $G$ is torsion-free, this theorem was independently proved in \cite{GM1,GM2}.

\begin{thm}\label{Fill}
Suppose that a group $G$ is hyperbolic relative to a collection of
subgroups $\Hl $. Then for most peripheral fillings of $G$, the
following holds.
\begin{enumerate}

\item[1)]  For each $\lambda\in \Lambda $, the natural map $H_\lambda\
/N_\lambda \to G(\N )$ is injective.

\item[2)] The quotient-group $G(\N )$ is hyperbolic relative to the
collection $\{H_\lambda /N_\lambda \} _{\lambda \in \Lambda }$.
\end{enumerate}
\end{thm}

The following statement plays a key role in our paper.

\begin{lem}\label{conj}
Let $G$ be a relatively hyperbolic group, $H$ -- a subgroup of $G$ and
$\alpha \in Aut(H)$. Suppose that there exists a loxodromic element
$g\in H$ such that $\alpha (g)$ is not conjugate to an element of
$E_G(g)$ in $G$. Then $\alpha $ does not preserve some normal
subgroup of $H$.
\end{lem}

\begin{proof}
Suppose that $G$ is hyperbolic relatively to $\Hl $. There are two
cases to consider.

{\bf Case 1.} Assume first that $\alpha (g) $ is loxodromic. Using
Lemma \ref{Eg} twice we obtain that $G$ is hyperbolic relatively to
$\Hl \cup \{ E_G(g), E_G(\alpha (g))\} $. Since $\langle g\rangle $
has finite index in $E_G(g)$, there is $m\ne 0$ such that $\langle
g^m\rangle $ (and each of its subgroups) is normal in $E_G(g)$. Let
$\mathcal F$ be the finite set provided by Theorem \ref{Fill} for
the peripheral system $\Hl \cup \{ E_G(g), E_G(\alpha(g))\} $.
Taking $p$ to be a sufficiently large multiple of $m$, we can ensure
the condition $\langle g^p\rangle \cap \mathcal F =\emptyset $. We
now consider the filling of $G$ with respect to the collection of
subgroups $\N $ consisting of the trivial subgroups of $H_\lambda$'s, the trivial subgroup of
$E_G(\alpha (g))$, and $\langle g^p\rangle\lhd E_G(g)$. By Theorem
\ref{Fill} elements $g$ and $\alpha (g)$ have orders $p$ and $\infty
$, respectively, in $Q=G/\ll g^p\rr^G$. Hence $\alpha $ does not
induce an automorphism on the natural image of $H$ in $Q$, i.e., it
does not preserve $\ll g^p\rr^G\cap H$.

{\bf Case 2.} Now suppose that $\alpha (g) $ is parabolic, i.e., it
is conjugate to an element of some peripheral subgroup $H_\lambda $.
Again, by Lemma \ref{Eg}, $G$ is hyperbolic relatively to $\Hl \cup \{
E_G(g)\} $. The rest of the proof is identical to that in Case 1.
The only difference is that Theorem \ref{Fill} is applied to the
collection of subgroups $\N $ consisting of trivial subgroups of
$H_\lambda$'s and $\langle g^p\rangle\lhd E_G(g)$ for some $p>0$.
\end{proof}

Theorem \ref{main} is a particular case of the following
result. (Recall that $\widehat H=H E_G(H)$.)

\begin{thm}\label{NormAutSubgr}
Let $G$ be a relatively hyperbolic group and let $H \le G$ be a
non-elementary subgroup such that $H^0\ne \emptyset $. Then for any
$\varphi\in Aut_n(H)$ there exists a map $\e : H \to E_G(H)$, whose
restriction to $\CH $ is trivial, and an element $w \in N_G (\widehat H
)$ such that for every $h \in H$, $\varphi(h)=w  h \e (h)
w^{-1}$.
\end{thm}

\begin{proof}
By Lemma \ref{conj}, $\varphi $ maps every loxodromic element $h\in H$ to a conjugate of an element of $E_G(h)$.
As $\langle h\rangle $ has finite index in $E_G(h)$, every element of infinite order in $E_G(h)$ is
commensurable with $h$ in $G$. In particular, $\varphi (h)\apg h$ for every $h\in H^0$.
Hence by Theorem \ref{thm:comm-aut} there is a map
$\e : H \to E_G(H)$, whose restriction to $\CH $ is a homomorphism,
and an element $w \in N_G (\widehat H )$ such that for every $h \in H$,
$\varphi(h)=w h \e (h) w^{-1}$. It remains to show
that $\e (h)=1$ for every $h\in \CH $.

By Proposition \ref{fi}, it suffices to show that $\e (h)=1$ for
all $h\in S_G(H)$. Suppose that $\varphi (h)=whrw^{-1}$ for some
$r\in E_G(H)\setminus \{ 1\}$. Take any integer $p \equiv 1~ ({\rm mod }\; |r|)$, where $|r|$ denotes the (finite) order of $r$ in $G$.
Note that $h$
commutes with $r$ as $h\in S_G(H)$. Thus $\varphi
(h^p)=wh^p rw^{-1}$. Since $\varphi$ should preserve $\ll h^p\rr ^G
\cap H$, we obtain $h^p r\in \ll h^p\rr ^G$. On the other hand,
$h^pr\in E_G(h)$. By Lemma \ref{Eg} we can join $E_G(h)$ to the
collection of the peripheral subgroups. Without loss of generality
we may assume that $p \gg 1$ so that the normal subgroup $N=\langle
h^p\rangle $ of $E_G(h)$ satisfies the requirement $N\cap \mathcal
F=\emptyset$ from Theorem \ref{Fill} (and Definition \ref{def:periph_fill}). Then by the first part of
Theorem \ref{Fill} we have $h^pr \in \ll h^p\rr ^G \cap E_G(h)
=\langle h^p\rangle $. Hence $r\in \langle h\rangle \cap E_G(H)=\{
1\}$, which contradicts $r\ne 1$.
\end{proof}

\begin{cor}\label{outn}
Let $H$ be a non-elementary subgroup of a relatively hyperbolic group
$G$ such that $H^0\ne \emptyset $. Then the following hold.
\begin{enumerate}
\item[a)] If $H$ has finite index in $N_G(HE_G(H))$, then $Out_n(H)$ is
finite.

\item[b)] If $H$ does not normalize any non-trivial finite subgroup of $G$, and $H=N_G(H)$, then
$Out _n(H)=\{ 1\}$.
\end{enumerate}
\end{cor}

\begin{proof} The argument is similar to the one used to prove Corollary \ref{cor:conj_aut}.
Observe that by Lemma \ref{EH}, $E_G(H)$ is a finite subgroup of $G$ normalized by $H$.
Therefore $H$ acts on $E_G(H)$ by conjugation,
and $C_H=C_H(E_G(H))$ has a finite index in $H$ as a kernel of this action.


Let $h_1, \dots , h_l$ be elements of $H$ such that $H=\bigsqcup_{i=1}^l C_H h_i$. By Theorem \ref{NormAutSubgr} we can argue
as in the proof of Corollary \ref{cor:conj_aut} to conclude that every
normal automorphism $\varphi$ of $H$ is uniquely determined by the images $\e(h_i)$ of $h_i$, $i=1, \dots,l$,
and by  the conjugating element $w \in N_G(\widehat H)$. As $E_G(H)$ is finite,
for each $i$ there are only finitely many possibilities for $\e(h_i)$, and since $|N_G(\widehat H):H|<\infty$, we can deduce that
$|Aut_n(H):Inn(H)|<\infty$.

Furthermore, if $H=N_G(H)$ and $H$ does not normalize any finite
normal subgroup of $G$, we obtain $E_G(H)=\{1\}$, $N_G(\widehat H)=N_G(H)=H$,
and $\CH =H$. Hence $Aut_n(H)=Inn (H)$ by Theorem
\ref{NormAutSubgr}. This completes the proof.
\end{proof}

The next lemma shows that Corollary \ref{cor1} holds for elementary groups.

\begin{lem}\label{vc}
Let $G$ be a virtually cyclic group. Then $Out(G)$ is finite.
\end{lem}

\begin{proof} If $G$ is finite the claim is trivial, so assume that $G$ is infinite.
Recall that every elementary group is ether finite-by-cyclic or finite-by-(infinite dihedral)
(see, for example,  \cite[Lemma 2.5]{FJ}). More precisely, as $G$ is infinite, the quotient $G/E(G)$
(where $E(G)$ is the maximal finite normal subgroup of $G$ given by Corollary \ref{KG}) is either
infinite cyclic or infinite dihedral. In both cases we have
\begin{equation} \label{eq:aut(g/eg)}\left|Aut(G/E(G)):Inn(G/E(G))\right|=2. \end{equation}

Every automorphism $\alpha \in Aut(G)$ induces an automorphism $\bar\alpha \in Aut(G/E(G))$.
This gives rise to a homomorphism $\xi: Aut(G)\to Aut(G/E(G))$. If $\alpha \in \ker(\xi)$, then for every
$x \in G$ there is $h=h(x) \in E(G)$ such that $\alpha(x)=x h$. By our assumptions, $G$ is generated by a finite set
of elements $\{x_i~|~i=1,\dots,n\}$ and  the automorphism
$\alpha$ is uniquely determined by the images $\alpha(x_i)$, $i=1,\dots,n$.
Since $|E(G)|<\infty$, for each $i$ there are only finitely many possibilities for $h(x_i)$. Therefore
the kernel of $\xi$ is finite. Evidently $\xi (Inn(G))=Inn (G/E(G))$, and by \eqref{eq:aut(g/eg)} we get
$|Aut(G):(Inn(G) \ker (\xi))| \le 2$ yielding that $|Out(G)|=|Aut(G):Inn(G)|<\infty$.
\end{proof}

\begin{proof}[Proof of Theorem \ref{main}]
Let us apply Theorem \ref{NormAutSubgr} to the case $G=H$. Then $E_G(H)=E(G)$, $\CH =C(G)$, $\widehat H=N_G(\widehat H)=G$, and the claim of Theorem \ref{main} follows immediately. 
\end{proof}

\begin{proof}[Proof of Corollary \ref{cor1}]
First, suppose that $G$ is elementary. In this case the first part of the corollary follows from Lemma \ref{vc}.
To derive the second claim of the corollary, we observe that since $G$ is non-cyclic and does not have non-trivial finite
normal subgroups, it must be infinite dihedral
(this follows from the structure of an elementary group -- see the proof of Lemma \ref{vc}).
Hence $G \cong \Z/2\Z * \Z/2\Z$ and, by Neshchadim's theorem \cite{Nesh},
$Out_n(G)=\{1\}$.

Thus we may assume that $G$ is non-elementary. In this case the corollary follows
from Theorem \ref{main} in the same way as Corollary \ref{outn} from Theorem \ref{NormAutSubgr}.
Alternatively it follows immediately from Corollary \ref{outn} applied to the case when $G=H$.
\end{proof}


\section{Free products and groups with infinitely many ends}


In order to prove Theorem \ref{fp} we need two more statements below.

\begin{lem}\label{noncom}
Assume that $G$ is a relatively hyperbolic group and $g, h$ are two
non-commensurable loxodromic elements. Then $g$ and $h$ are
non-commensurable and loxodromic in most peripheral fillings of $G$.
\end{lem}

\begin{proof}
Suppose that $G$ is hyperbolic relative to a collection of subgroups
$\Hl $. Applying Lemma \ref{Eg} twice we obtain that $G$ is
hyperbolic relative to the new collection $\Hl \cup \{ E_1, E_2\}$, where $E_1=E_G(g)$,
$E_2=E_G(h)$. Let $\mathcal F_1$ and $\mathcal F_2$ be the finite
subsets provided by Theorem \ref{Fill} for the collections of
peripheral subgroups $\Hl$ and $\Hl \cup \{ E_1, E_2\}$,
respectively. Set $\mathcal F=\mathcal F_1\cup \mathcal F_2$.

Consider any collection of subgroups $N_\lambda \lhd H_\lambda $
such that $N_\lambda \cap \mathcal F=\emptyset$, $\lambda \in \Lambda$. By Theorem \ref{Fill}, the filling of $G$
with respect to the collection of normal subgroups $\N$,
consisting of $N_\lambda \lhd H_\lambda $ for $\lambda \in \Lambda $
and the trivial subgroups of $E_1$, $E_2$, is hyperbolic relative to
$\{ H_\lambda /N_\lambda \}_{\lambda \in \Lambda }\cup \{ E_1,
E_2\}$ as well as relative to $\{ H_\lambda /N_\lambda \}_{\lambda
\in \Lambda }$. (We keep the same notation for the isomorphic images
of $E_1, E_2$ in $G(\N)$ and the elements $g, h$.)

In particular, $E_1\cap E_2^t$ is finite for every $t\in
G(\N)$. Clearly this implies that $g$ and $h$ are not
commensurable in $G$. Similarly $g$ and $h$ are not conjugate to any
elements of the subgroups $H_\lambda /N_\lambda $, $\lambda \in
\Lambda $, of $G(\N)$. Thus $g$ and $h$ are loxodromic in $G(\N)$ with
respect to the peripheral collection $\{ H_\lambda /N_\lambda
\}_{\lambda \in \Lambda }$. As $\mathcal F$ is finite, $g$ and $h$
are non-commensurable and loxodromic in most peripheral fillings of
$G$ (with respect to the peripheral structure $\Hl $).
\end{proof}

The proof of Theorem \ref{fp} uses the following lemma, which is an immediate corollary of \cite[Lemma 3]{Yed}.
(Recall that the {\it Cartesian subgroup} of a free product $A*B$ is, by definition, the kernel of the natural epimorphism
$A*B \to A \times B$.)

\begin{lem}\label{Yed} Let $G=A\ast B$, where $A$ and $B$ are finite groups. Let
$u$, $v$ be non-commensurable elements of the Cartesian subgroup $C$ of $G$. Suppose that $u=a^k$, $v=b^l$
for some positive integers $k,l$, where $a$, $b$ are not proper powers. Assume also that $a^k$
(respectively, $b^l$) is the smallest non-zero power of $a$ (respectively, $b$) that belongs to $C$.
Then there exists a finite quotient-group $Q$ of $G$ such that the images of $u$ and $v$ have different orders in $Q$.
\end{lem}

\begin{proof}[Proof of Theorem \ref{fp}]
Let $G$ be a non-trivial free product, i.e., $G=A\ast B$, where both
$A$ and $B$ are non-trivial. Then $G$ is hyperbolic with
respect to $\{A,B\}$ (the finite sets $X$ and $\mathcal{R}$, from the definition of relative
hyperbolicity in Section \ref{sec:prelim}, can be taken to be empty; the isoperimetric constant $C$ for the corresponding relative
presentation of $G$ will then be equal to zero).  In what follows, we will fix this as a system of peripheral
subgroups of $G$.

If $|A|=|B|=2$, the proof is an easy exercise. It also follows from the main result of \cite{Nesh},
stating that every normal automorphism of a non-trivial free product is inner,
and the observation that every non-trivial normal subgroup of the infinite dihedral group is
of finite index.

Thus we may assume that $G$ is non-elementary. Suppose that there exists an automorphism $\alpha \in Aut_n^f(G)\setminus Inn(G)$.
Note that $E(G)=\{ 1\} $ because $G$, as a non-trivial free product, cannot contain non-trivial finite normal subgroups.
Since $\alpha$ is not an inner automorphism of $G$, it follows from Corollary \ref{cor:descr_comm_aut} that $\alpha $ is not commensurating.
Therefore, by Corollary \ref{cor:comm_aut_def} and Lemma \ref{lem:lox_non-comm} (applied to the case when $H=G$),
there is a loxodromic element $g\in G$ such that  $h=\alpha (g)$ is also loxodromic
and is not commensurable with $g$. Further, by Lemma \ref{noncom} there exist finite index normal
subgroups $M\lhd A$ and $N\lhd B$ such that the natural images
$\bar g$, $\bar h$  of $g$ and $h$, respectively, are not commensurable in
$\overline{G}=A/M\ast B/N$. Without loss of generality we may assume that $\overline{G}$ is non-elementary.

Since $\overline{G}$ is a
free product of two finite groups, it is residually finite. Therefore the kernel $K$ of the natural homomorphism $G\to \overline{G}$ is an intersection of finite
index normal subgroups of $G$.  As $\alpha \in Aut_n^f(G)$, $\alpha
$ stabilizes $K$. Hence $\alpha $ induces an automorphism
$\bar\alpha $ of $\overline{G}$.

Let $\bar g=a^k$, where $k$ is a positive integer and $a$ is not a proper power. Clearly $b=\bar\alpha (a)$ is
not a proper power as well and $b^k=\bar h$. Evidently $b^p=\bar\alpha (a^p)$ is not commensurable to $a^p$
for any non-zero integer $p$. Let $C$ denote the Cartesian subgroup
of $\overline{G}$. Then $|\overline{G}:C|<\infty$, and replacing $\bar g$ with another positive power of $a$, if necessary,
we may assume that $k>0$ and $\bar g=a^k$ is the smallest non-zero power of $a$ that belongs to $C$.
Again, since $|\overline{G}:C|<\infty $, $\bar \alpha $ preserves $C$.
In particular, $\bar h=b^k$ is the smallest power of $b$ that belongs to $C$.

By Lemma \ref{Yed} there exists a finite index normal subgroup $K$
of $\overline{G}$ such that the images of $\bar g$ and $\bar h$ have
different orders in $\overline{G}/K$. Therefore $\bar\alpha $ does
not induce an automorphism on $\overline{G}/K$. Obviously this means that $\alpha $ does not
preserve the full preimage of $K$ in $G$, which contradicts our
assumption that $\alpha \in Aut_n^f(G)$.
\end{proof}

The following lemma is well known and is easy to prove (see, for
example, \cite[Lemma 5.4]{GL}).

\begin{lem}\label{outnormsub}
Suppose that $G$ is a finitely generated group and $N$ is a centerless normal
subgroup of finite index in $G$. Then some finite index subgroup of
$Out (G)$ is isomorphic to a quotient of a subgroup of $Out (N)$ by
a finite normal subgroup. In particular, if $Out(N)$ is residually finite,
then $Out (G)$ is residually finite.
\end{lem}

The next observation is trivial.

\begin{lem}\label{finorb}
Suppose that a group $G$ acts on a set $\mathcal M$ faithfully with finite orbits. Then $G$ is residually finite.
\end{lem}

\begin{proof}
Given $g\in G$, let $s\in \mathcal M$ be an element such that $g(s)\ne s$. Then the natural map from $G$
to the symmetric group on the orbit of $s$ provides us with a finite quotient of $G$, where the image of $g$ is non-trivial.
\end{proof}

\begin{proof}[Proof of Theorem \ref{infends}]
Since the outer automorphism group of any virtually cyclic group is finite (see Lemma \ref{vc}),
we can assume that $G$ has infinitely many ends.

By Stallings's Theorem (\cite{Stall71,Stall68}) there is a finite group $S$ such that $G$ splits as an amalgamated free product
$A\ast_S B$ or an $HNN$-extension $A\ast _S$, where
$(|A:S|-1)(|B:S|-1)\ge 2$ in the first case and $|A:S_i|\ge 2$, $i=1,2$, in the
second case (where $S_1$ and $S_2$ are the two  associated isomorphic copies of $S$ in $A$). 
Since $G$ is residually finite and $S$ is finite, there exists a finite
index normal subgroup $N\lhd G$ such that $N\cap S=\{1\}$ if $G=A\ast_S B$, or $N \cap S_i = \{1\}$ for $i=1,2$, if $G=A\ast _S$. Note that
the quotient of the Bass-Serre tree for $G$ modulo the action of $N$
is finite and the edge stabilizers in $N$ are trivial. The
Bass-Serre structure theorem for groups acting on trees (see \cite{Serre}) yields
a splitting of $N$ into a non-trivial free product. In particular, $N$ is
centerless.

The group $Aut (N)$ naturally acts on the set $\mathcal M$ of finite index normal
subgroups of $N$ and $Aut_n^f(N)$ is the kernel of this action. By Theorem
\ref{fp}, $Aut_n^f(N)=Inn(N)$. Therefore, $Aut(N)/Aut_n^f(N)=Aut(N)/Inn(N)=Out(N)$ acts on $\mathcal M$ faithfully.
Since $N$ is finitely generated, there are only finitely many subgroups of a given finite index in $N$,
thus all orbits of the action of $Out(N)$ on $\mathcal M$ are finite. Hence $Out(N)$ is residually finite by Lemma \ref{finorb}.
The claim of the theorem is now a consequence of Lemma \ref{outnormsub}.
\end{proof}

\end{document}